\documentclass[11pt]{amsart}

\usepackage{amsmath,mathtools}
\usepackage[dvipsnames,svgnames]{xcolor}
\usepackage{amsmath, amsthm, amssymb, amsfonts, enumerate}
\usepackage[colorlinks=true,linkcolor=blue,urlcolor=blue]{hyperref}
\usepackage{dsfont}
\usepackage{color}
\usepackage{geometry}
\usepackage [latin1]{inputenc}
\usepackage{comment}
\usepackage{subcaption}
\usepackage{epstopdf}
\usepackage{bbm}
\usepackage{geometry}
\usepackage{setspace}

\usepackage{amsfonts}
\usepackage{amsfonts}
\usepackage{textcomp}

\usepackage{amssymb}
\usepackage{float}
\usepackage{tikz}
\usepackage{epsfig}
\usepackage{amsmath}
\usepackage[english]{babel}
\usepackage{a4}
\usepackage{enumerate}
\usepackage{amsmath}

\newcommand{\E}{\textsf{E}}
\renewcommand{\P}{\textsf{P}}
\newcommand{\one}{\text{$\mathbbm{1}$}}
\renewcommand{\hat}{\widehat}
\renewcommand{\tilde}{\widetilde}

\renewcommand{\epsilon}{\varepsilon}

\newtheorem{theorem}{Theorem}[section]
\newtheorem{remark}[theorem]{Remark}
\newtheorem{assumption}[theorem]{Assumption}
\newtheorem{lemma}[theorem]{Lemma}

\newtheorem{corollary}[theorem]{Corollary}

\newtheorem{definition}[theorem]{Definition}

\numberwithin{equation}{section}

\title[A Stochastic Non-zero-sum game of controlling the debt-to-GDP ratio]{{A Stochastic Non-zero-sum game of controlling the debt-to-GDP ratio}}

\author[Dammann]{Felix Dammann$^{\ddag} $}
\thanks{$^\ddag$Center for Mathematical Economics (IMW), Bielefeld University, Universit\"{a}tsstrasse 25, 33615, Bielefeld, Germany. The author gratefully acknowledges funding by the Deutsche Forschungsgemeinschaft (DFG, German Research Foundation) -- Project-ID 317210226 -- SFB 1283.}
\author[Rodosthenous]{Neofytos Rodosthenous$^*$}
\thanks{$^*$Department of Mathematics, University College London, Gower St, London WC1E 6BT, UK}% Email: n.rodosthenous@ucl.ac.uk}
\author[Villeneuve]{St\'{e}phane Villeneuve$^\dagger$}
\thanks{$^\dagger$Toulouse School of Economics, University of Toulouse  Capitole, 1 esplanade de l'universit\'{e}, 31000 Toulouse, France, the authors acknowledge funding from  ANR under grant ANR-17-EUR-0010 (Investissements d'Avenir program) and gratefully thanks the FdR-SCOR
{\it Chaire March\'e des risques et cr\'eation de valeurs}. }

\date{\today}

\begin{document}
\begin{abstract} 
We introduce a non-zero-sum game between a government and a legislative body to study the optimal level of debt.  
Each player, with different time preferences, can intervene on the stochastic dynamics of the debt-to-GDP ratio via singular stochastic controls, in view of minimizing non-continuously differentiable running costs.
We completely  characterise Nash equilibria in the class of Skorokhod-reflection-type policies. 
We highlight the importance of different time preferences resulting in qualitatively different type of equilibria. 
In particular, we show that, 
 while it is always optimal for the government to devise an appropriate debt issuance policy,
the legislator should optimally impose a debt ceiling only under relatively low discount rates and a laissez-faire policy can be optimal for high values of the legislator's discount rate.
\end{abstract}

\maketitle

\text{\it Keywords}: non-zero-sum game, singular stochastic control, free-boundary problem, debt-to-GDP ratio.

\smallskip

{\text{\it MSC2010 subject classification}}: 91A05, 93E20, 90B05

{\text{\it JEL subject classification}}: C61, C73

\section{Introduction}

There is probably no more topical issue in macroeconomics than the determination of the optimal level of debt that favours both its sustainability and the long-term growth of an economy. 
Although paramount and extensively studied in the literature (see Barro \cite{Ba} and Dornbusch and Draghi \cite{DoDr} for a general presentation), the question of the optimal debt level has not yet received clear theoretical foundations. 
This lack of a consensual theoretical framework has led to the implementation of exogenous mechanisms to monitor the level of debt. 
In the USA, one of these mechanisms is the statutory debt ceiling which restricts the amount of debt a government can be permitted to issue\footnote{Within the European Community, a similar mechanism exists since the Maastricht treaty set $60\%$ as the upper bound for the debt-to-GDP ratio for members of the European Union.}. \\
The traditional analysis of public debt has shown that high public debt has a negative effect on long-term economic growth, thus giving an argument to debt ceiling advocates. 
Indeed, a high level of debt generates high risk premiums that reflect creditors' doubts about the government's ability to refinance itself. 
Being unable not only to repay its debts but also to pay for the excess of its expenditures over its revenues, the government must then immediately balance its budget by taking exceptional measures, like increasing taxes and/or cutting its investments, which can have a dramatic impact on growth. 

Motivated by this, the theoretical literature  on debt management problems has focused on the stochastic control problem faced by a single decision maker to determine the optimal debt reduction policy and thus the debt ceiling, i.e.~the level of debt-to-GDP ratio (also called ``debt ratio") at which the government should intervene in order to reduce it. 
In Cadenillas and Huam\'an-Aguilar \cite{CaHu, CaHu2}, the debt ratio evolves as a controlled one-dimensional geometric Brownian motion that can be reduced via singular and bounded variation controls, respectively, in order to minimize the expected total costs resulting from the instantaneous cost of the debt ratio and intervention costs. 
Ferrari \cite{Fe} studies the optimal debt ratio reduction problem posed as a fully two-dimensional singular stochastic control problem, where the government takes into consideration the evolution of the inflation rate (evolving as an uncontrolled diffusion process) of the country as well. 
Callegaro et al.~\cite{CaCeFe} consider a model with partial observation, where the growth rate of GDP follows an unobservable Markov chain.

On the other hand, the positive effect of a high level of public debt on growth should not be overlooked, since public investments in social policies, education, healthcare, justice, research, and infrastructure help private initiatives to develop effectively. 
As Blanchard observed in his presidential address to the American Economic Association \cite{Bl}, as long as the interest rate is lower than the growth rate, a large deficit can be allowed without decreasing the debt ratio.
Notable contributions that consider this ambivalent effect of public debt include Ferrari and Rodosthenous \cite{FeRo}, who model the growth rate of GDP by a continuous-time Markov chain and the government is allowed to both decrease and increase the debt ratio, as well as Brachetta and Ceci \cite{BrCe}, who consider a model of regular controls where interventions via fiscal policies affect the public debt and the GDP growth rate at the same time. 

A common feature in all aforementioned papers is that an optimal level of debt ceiling is endogenously determined. 
However, the potential political game between a government and a legislative body (e.g.~the US Congress), whose political interests may be divergent, has not been considered so far. 
In our paper, we contribute to the literature on debt management by proposing a game that incorporates the opposing interests and strategic interaction between these two players. %, and thus allows us to model the ambivalent effect of public debt studied in the literature. 
Our modeling framework results in a non-zero-sum game between a government, whose mandate is to manage its public debt issuance policy to finance its spending, 
while a legislative body is concerned with imposing a mechanism limiting the amount of debt to avoid a potential debt crisis. 
In mathematical terms, each player exerts a monotone control to set the path of the stochastically evolving debt ratio. 
The first player (government) can increase the level of debt ratio by exerting its control, while the
second player (legislative body) can decrease the level of debt ratio by implementing exceptional measures. 
Each player aims at minimizing their own total cost functional and we allow the rate of increase/reduction of each player to be unbounded and have an instantaneous effect on the debt ratio. 
Consequently, this leads to the formulation of a stochastic non-zero-sum game of singular controls.

Even though there exists a considerable literature on one-player singular control problems (e.g. Bather and Chernoff \cite{BaCh}, Benes et al.~\cite{Be}, Karatzas \cite{Ka} and many others), the literature on non-zero-sum stochastic games with singular controls is still limited. 
Kwon and Zhang \cite{Kw} study a game of competitive market share control, where each player can make irreversible investment decisions via singular controls as well as decide to exit the market, and obtain and characterise Markov perfect equilibria. 
Our work is more closely related to De Angelis and Ferrari \cite{DeAFe2}, who prove the existence of a Nash equilibrium in the class of Skorokhod-reflection policies, by establishing a new connection of a non-zero-sum game of monotone controls with a non-zero-sum stopping game. 
In order to achieve this connection, it is necessary that both players have the same discount factor (or equivalently time preferences) and that the running cost is a differentiable function -- same assumptions are imposed also in \cite{Kw}. 
In our paper, we relax both of these assumptions, by considering different time preferences for each player and a non-differentiable running cost function. 
This results in the need for a different methodology to prove the existence of a Nash equilibrium.

We first study separately two coupled constrained stochastic control problems faced by the two players, and then search for Nash equilibria in the game, i.e.~where no player has an incentive to deviate unilaterally. 
From the government's perspective, assuming that the legislative body imposes a debt ceiling $b$ (or equivalently a Skorokhod reflection policy for the debt ratio process at $b$), we investigate the optimal policy of the government for issuing new public debt. 
To this end, we establish a connection between this constrained stochastic control problem and a free-boundary problem, that we solve via a guess-and-verify approach. 
As a result, we show that the best debt issuance policy is to reflect the debt ratio process upwards at a level $a(b)$, which clearly depends on the imposed debt ceiling mechanism by the legislator. 
Consequently, considering that the government is going to use a debt issuance policy of reflecting the debt ratio process at a level $a$, the legislator should decide on whether to impose a debt ceiling. In particular, if there is already a statutory exogenous debt ceiling, is it optimal to raise it and by how much? 
Our results suggest that a debt ceiling $b(a)$ should indeed be imposed for a specific range of the legislative body's time preference rates. For larger values -- implying that the legislator discounts future costs more heavily -- the legislator's optimal strategy is in fact realised by a laissez-faire policy in which no debt ceiling mechanism is imposed.

The main contribution of our paper is to eventually prove the existence and uniqueness of a Nash equilibrium in the game. 
We show that -- depending on the legislator's time preference rate $\lambda$ -- two qualitatively different Nash equilibria exist in the game. More precisely, for specific values of $\lambda$ the Nash equilibrium prescribes that the debt ratio is kept inside the interval $[a^*,b^*]$ with the minimal cost, associated to Skorokhod reflection policies. 
On the other hand, for large values of the legislative body's time preference rate, the legislative body should optimally not intervene, and an associated Nash equilibrium without debt ceiling is proved to hold. Interestingly, we prove that the optimality of adopting such a strategy relies solely on the legislative body's time preferences compared to the parameter constellation in the model -- it does {\it not depend} on the actions of the opposing player -- see specifically our results in Section \ref{Sec: ECB Case I}. 
When the discount rate of future costs is high, the consideration of the risk of a debt crisis in the future is too low to make the implementation of a debt ceiling mechanism optimal.
%We thus split our search for Nash equilibria in based on the magnitude of the legislative body's time preferences.

The paper is organized as follows. 
Section 2 describes the setting and the two problems faced by the two players. 
In Section 3, we solve the government constrained control problem by distinguishing two cases: the legislative body does not intervene or forces the government to keep its debt ratio below a debt ceiling %fixed leve 
$b>0$. 
In both cases, we show that the government can devise an optimal debt issuance policy when the debt ratio is sufficiently low.
In Section 4, we solve the constrained control problem of the legislative body and find its best response strategy to the above governmental policy.
We show that the magnitude of the legislative body's time preference rate plays a crucial role, which is a key result of our analysis. 
In particular, if it is relatively small, a debt ceiling mechanism is optimal, while if it is relatively large, the legislative body should not set a debt ceiling at all.
We prove the existence and uniqueness of a Nash Equilibrium in the class of Skorokhod-reflection policies in Section 5. 
Finally, Section 6 is devoted to a comparative statics analysis; we explore how the optimal debt issuance policy and debt ceiling mechanism are affected by changes in the model parameters. 
In particular, we are able to quantify the transition between a legislative body's optimal intervention and non-intervention regimes.

\section{Setting and Problem Formulation}\label{Sec: Setting and Problem Formulation}
\subsection{Motivation for the model} The model applies to a government that has to finance its expenditures through public debt, under the control of a legislative body. The nominal debt grows at rate $r$, i.e. it evolves according to
\begin{align*}
d D_t = r D_t dt, 
\qquad t \geq 0, 
\end{align*}
%for times $t \geq 0$, 
in the absence of any intervention, where we denote by $r \in \mathbb{R}$ the interest rate on government debt. 
When designing its economic policy, the government can choose to increase
the current level of the debt by a new issuance.
Denoting by $\xi_t$ the
cumulative percentage of debt increase by the government up to time $t \geq 0$, the dynamics of the adjusted debt reads as 
\begin{align*}
d D_t = r D_t dt + D_t  \,\circ_u\,d \xi_t,
\qquad t \geq 0, 
\end{align*}
where the latter integral %(and others of the same type in the followingmotivation) 
will be defined later in Section \ref{probform}.
The GDP follows the stochastic exogenous dynamics
\begin{align*}
d G_t = g G_t dt + \sigma G_t d \hat{W}_t, 
\qquad t \geq 0, 
\end{align*}
%for $t \geq 0$, 
in the absence of any intervention, where we denote by $g \in \mathbb{R}$ the growth rate of the GDP and by $\hat{W}$ a standard one-dimensional Brownian motion. 
A legislative body %International financial institutions 
can implement a liberalisation policy in order to boost the GDP by forcing the government to favor the job market, moderating social insurance programs, reducing burdensome regulations, lowering the marginal tax rate and privatizing businesses.
Denoting by $\eta_t$ the cumulative percentage of GDP increase implemented by the legislative body up to time $t \geq 0$, the dynamics of the adjusted GDP read as 
\begin{align*}
d G_t = g G_t dt + \sigma G_t d \hat{W}_t + G_t  \,\circ_u\, d \eta_t,
\qquad t \geq 0. 
\end{align*}
Hence, we may conclude that the dynamics of the debt-to-GDP ratio, obtained via the use of
It\^{o}'s formula on $X := D/G$, evolves according to 
\begin{align*}
d X_t = (r- g) X_t dt + \sigma X_t (\sigma dt - d\hat{W}_t) + X_t  \,\circ_u\, d \xi_t - X_t  \,\circ_d\, d \eta_t,
\qquad t \geq 0, 
\end{align*}
where the latter integral will also be defined later in Section \ref{probform}.
Without loss of generality, a change of measure to one under which $W_t := \sigma t - \hat{W}_t$ is a Brownian
motion, will allow for the following problem formulation.

\subsection{Problem formulation} 
\label{probform}
Let $(\Omega, \mathcal{F}, \P)$ be a complete probability space accommodating a one-dimensional Brownian motion $W := (W_t)_{t\geq 0}$. We denote by $\mathbb{F} := \{ \mathcal{F}_t , t \geq 0 \}$ the filtration generated by $W$ augmented by $\P$-null sets. In absence of any interventions, the debt-to-GDP ratio (also called ``debt ratio'') evolves according to the stochastic differential equation (SDE) 
\begin{align}\label{Dyn: DRatio Uncontrolled}
d X_t^0 = (r - g)X_t^0 dt + \sigma X_t^0 d W_t , \qquad X_t^0 = x >0.
\end{align} 
The classical macroeconomic dynamics of the debt ratio, see e.g.~\cite{BlFi}, are simply the deterministic version of \eqref{Dyn: DRatio Uncontrolled} with $\sigma = 0$. 
When increasing the current debt ratio level by $\epsilon > 0$ percentage points, the debt ratio exhibits a jump 
\begin{align*}
\Delta X_t = X_t - X_{t-} = \epsilon X_{t-} .
\end{align*}
Hence, for small $\epsilon>0$, we can associate a governmental intervention on the debt ratio with $X_t = (1 + \epsilon)X_{t-} \approx e^\epsilon X_{t-}$. Furthermore, interpreting an intervention $\Delta \xi_t$ as a sequence of $N$ individual interventions of size $\epsilon = \Delta \xi_t / N$ we have $X_t = e^{N \epsilon} X_{t-} = e^{\Delta \xi_t} X_{t-}$, for $N$ large enough. 
We can thus model the controlled debt ratio dynamics (by arguing similarly for the interventions of the legislative body) via 
\begin{align}\label{Dyn: DRatio Controlled}
d X_t^{\xi,\eta} = (r-g) X_t^{\xi,\eta} dt + \sigma X_t^{\xi,\eta} d W_t + X_t^{\xi,\eta} \circ_u d \xi_t - X_t^{\xi,\eta} \circ_d d \eta_t,
\qquad t \geq 0, 
\end{align}
where the operators $\circ_u$ and $\circ_d$ are defined as 
(see also \cite{DeAFe2})
\begin{align}\label{Def: Operator circ up}
\begin{split}
X_t^{\xi,\eta} \circ_u d \xi_t &= X_t^{\xi,\eta} d \xi_t^c + X_{t-}^{\xi,\eta} \int_0^{\Delta \xi_t} e^u du =  X_t^{\xi,\eta} d \xi_t^c + X_{t-}^{\xi,\eta} \big[ e^{\Delta \xi_t} - 1 \big], \\ 
%\end{align}
%and 
%\begin{align}\label{Def: Operator circ down}
X_t^{\xi,\eta} \circ_d d \eta_t &= X_t^{\xi,\eta} d \eta_t^c + X_{t-}^{\xi,\eta} \int_0^{\Delta \xi_t} e^{-u} du =  X_t^{\xi,\eta} d \eta_t^c + X_{t-}^{\xi,\eta} \big[ 1 - e^{\Delta \xi_t}  \big].
\end{split}
\end{align}
Here, $\xi^c$ (resp.,~$\eta^c$) denotes the continuous part of the process $\xi$ (resp.,~$\eta$). 

Using It\^{o}'s formula we can verify that the solution to \eqref{Dyn: DRatio Controlled} starting at time zero from level $x > 0$ is given by 
\begin{align} \label{Xxieta}
X_t^{\xi,\eta} = x \exp \Big( \Big( r - g- \frac{1}{2} \sigma^2 \Big) t + \sigma W_t + \xi_t - \eta_t \Big) = X_t^0 \exp(\xi_t - \eta_t),
\qquad t \geq 0, 
\end{align}
where $X_t^0$ denotes the solution to \eqref{Dyn: DRatio Uncontrolled}. 
Notice that the impact of interventions by the government and legislative body are of multiplicative structure and additive to the logarithm of the debt ratio. 
Finally, we introduce the process $Y^{\xi,\eta}_t:=\ln{X_t^{\xi,\eta}}$, such that $Y_t^0 = \ln{X^0_t}$ from \eqref{Dyn: DRatio Uncontrolled}, to obtain from \eqref{Xxieta} that
\begin{align}\label{logX}
Y^{\xi,\eta}_t = Y_t^0 + \xi_t - \eta_t,
\qquad t \geq 0. 
\end{align}

In accordance with our reasoning above, $\xi_t$ denotes the cumulative percentage amount of debt increase by the government and $\eta_t$ denotes the cumulative percentage amount of debt decrease by the legislative body, up to time $t \geq 0$. 
It is therefore natural to model them as nondecreasing stochastic processes, adapted with respect to the available flow of information $\mathbb{F}$. Hence we take $\xi$ and $\eta$ in the set 
\begin{align*}
\mathcal{U} := \Big\{ \upsilon: \Omega \times \mathbb{R}_+ \to \mathbb{R}_+ :~  ~(\upsilon_t)_{t \geq 0} ~\mathbb{F}\text{-adapted, nondecreasing, c\`{a}dl\`{a}g, and } \upsilon_{0-} = 0 \Big\}.
\end{align*}

\subsection*{The problem of the government}
In this framework, the government is facing a potential debt ceiling (or debt limit) as a hard constraint imposed by a legislative body, when the country's debt ratio is too high. 
In other words, the government has an exogenous factor, namely a debt ratio ceiling $b$, to take
into consideration when designing its economic policy. This is the level at which a legislative body will demand the decrease of the debt ratio and the
adoption of liberalisation policies by the government. 
In the following, we assume that having a debt level $X_t^{\xi,\eta}$ at time $t \geq 0$, the government incurs an
instantaneous cost $h(X_t^{\xi,\eta})$. This may be interpreted as an opportunity cost resulting from
having less room for financing public investments. We make the following standing assumption.

\begin{assumption}\label{Assumption: Cost function h} 
The instantaneous (running) cost function $h: \mathbb{R}_+ \to \mathbb{R}_+$ satisfies: 
\begin{itemize}
\item[(i)] $x \mapsto h(x)$ is strictly convex, continuously differentiable and increasing on $[0,\infty)$; 
\item[(ii)] the derivative $h'$ of $h$ satisfies $\lim_{x \to 0}h'(x) = 0$ and $\lim_{x \to \infty} h'(x) = + \infty$;
\item[(iii)] there exists $p > 1$, $K_1 > 0$ such that 
\begin{align*}
h(x) \leq K_1 (1 + \vert x \vert^p ) , \quad x \in \mathbb{R}.
\end{align*}
\end{itemize}
\end{assumption}
\begin{remark} \label{remh}
It is worth noticing that a cost function of the form $h(x) = \frac12 x^2$ for $x > 0$ satisfies Assumption \ref{Assumption: Cost function h}. 
Notice that $h(0) = 0$ together with $h'(0) = 0$ imply that any infinitesimal amount of debt does not generate holding costs for the country;
indeed, $h(\epsilon) \approx h'(0)\epsilon = 0$. 
If one wishes to obtain closed-form solutions, a specific function $h$ must be chosen according to Assumption \ref{Assumption: Cost function h}; our choice will be precisely the above one.
\end{remark}
Moreover, whenever a legislative body decides to impose a debt ceiling mechanism, the government incurs a
proportional cost to the amount of debt reduction (see also \cite{CaHu}, \cite{Fe} and \cite{FeRo}). This might be seen as a measure of the social and financial consequences, or repercussions for the financial stability
of households and individuals, deriving from the enforcement of debt-reduction policies. The
associated constant marginal cost $c_1 > 0$ allows to express it in monetary terms.
Finally, the government's main aim is to increase the current level of debt ratio through
public investments, e.g. investments in infrastructure, healthcare, education and research, etc.
We assume that this has a positive political, social and financial effect, thus overall reduces the total
expected ``costs" of the government. The marginal benefit of increasing the debt ratio is a
strictly positive constant $c_2 > 0$.
From the point of view of the government, assuming that it discounts at a rate $\rho > 0$, the
total expected cost functional, net of investment benefits, is thus given by
\begin{align}\label{Objective: Government}
\mathcal{J}_{x,\eta}(\xi) := \E_x \Big[ \int_0^\infty e^{- \rho t} h(X_t^{\xi,\eta})dt + c_1 \int_0^\infty e^{- \rho t} X_t^{\xi,\eta} \circ_d d \eta_t - c_2 \int_0^\infty e^{- \rho t}X_t^{\xi,\eta} \circ_u d \xi_t \Big],
\end{align}
where, for any $x \in \mathbb{R}_+$, $\E_x$ denotes the expectation under the measure $\P_x ( \cdot) := \P(\cdot \mid X^{\xi,\eta}_{0_{-}} = x)$.

\subsection*{The problem of the legislative body}
On one hand, the legislative body (e.g. Congress) would like governments
to ideally keep their country's debt ratio at low levels to maintain a low probability
of default and a feasible borrowing from the markets. Even though countries that can print
their own currency cannot default on their debts, there are many countries that do not control
their own monetary policy, e.g. EU members who rely on the European Central Bank (ECB), or countries that hold
large amounts of foreign denominated debts, e.g.~ Argentina (who defaulted on US government bonds). Several levels $m>0$ defining the ``healthy" region $[0,m]$ of relatively ``low" debt ratio have
been used in the last decades, e.g. $m = 60\%$ is the Maastricht Treaty's reference value of 1992
for all EU countries, or $m = 77\%$ is the threshold found by researchers at the World Bank \cite{CaGrKo} for developed economies and $m = 64\%$ for emerging markets.\footnote{This study goes a step further to quantify the economic cost per percentage point the debt ratio exceeds $m$ (see also \cite{ReReRo}  for an empirical study on the effect of high debt towards private investments' crowding out and a low subsequent growth)}

When the debt ratio $X^{\xi,\eta}$ exceeds this pre-specified value $m>0$, the legislative body would face social and political pressure, which may lead to the implementation of liberalisation policies in order to decrease the level of $X^{\xi,\eta}$ via a control strategy $\eta$. 
This could, for example, be done by setting a debt ceiling $b$. 
This debt ceiling $b$ is expected to be bigger than $m$, since imposing structural adjustment programs on countries or restricting further borrowing by governments, is costly for the legislative body and the associated marginal cost is $\kappa > 0$. 
From the point of view of such a legislative body, assuming that it discounts at a rate $\lambda > 0$, we model the expected cost functional as\footnote{We again highlight the fact that the legislative body discounts with a different discount rate than the government, which can be interpreted as different time preferences. 
Moreover, the running cost function inside the first integral is non-differentiable, which does not satisfy the assumptions in \cite{DeAFe2}, thus the link between non-zero-sum games of singular controls and optimal stopping, developed therein, breaks down.} %observed in \cite{DeAFe2}.}
\begin{align}\label{Objective: ECB New}
\mathcal{I}_{x,\xi} (\eta) := \E_x \Big[ \int_0^\infty e^{-\lambda t} \alpha (X_t^{\xi ,\eta} - m)^+dt + \kappa \int_0^\infty e^{- \lambda t} X_t^{\xi,\eta} \circ_d d \eta_t \Big]. %\\
%\mathcal{I}_{x,\xi}^+ (\eta) := \E_x \Big[ \int_0^\infty e^{-\lambda t} \alpha (X_t^{\xi ,\eta} - m)^+ dt + \kappa \int_0^\infty e^{- \lambda t}  X_t^{\xi,\eta} \circ_d d \eta \Big], 
\end{align}

When a government wants to reduce its public deficit, it has, in simple terms, a choice between increasing tax revenues while keeping expenditures constant, or reducing public expenditures with stable tax revenues. The second choice is usually the more difficult to make: public spending is sometimes structural (for example, the payment of civil servants' salaries) and therefore incompressible in the short term. This is why, when seeking to reduce public deficits, one most frequently turns to taxation. Hence, $\alpha>0$ can be interpreted as a country tax compliance factor, the smaller $\alpha$ is the bigger is the willingness to pay tax, if needed in the future. When this factor is low as it is in Denmark for instance, the legislative body has thus less social pressure to reduce the debt ratio.

\subsection*{Debt ceiling mechanism as a non-zero-sum game of singular controls.}
In our analysis, we restrict our attention to controls producing finite payoffs, which includes the realistic assumption that both players will not use an economic policy leading to infinite cost and/or benefit of interventions. 
Moreover, we note that the definition of the integrals with respect to the controls, as specified in \eqref{Def: Operator circ up}, %-\eqref{Def: Operator circ down} 
requires some attention since simultaneous jumps of $\xi$ and $\eta$ may be difficult to handle.

Given that the debt ratio is always a positive number, we therefore consider pairs $(\xi,\eta) \in \mathcal{U} \times \mathcal{U}$ such that 
\begin{align}\label{Equ: Obj 1 finite}
\begin{split} 
&\E_x \Big[ \int_0^\infty e^{- \rho t} h(X_t^{\xi,\eta})dt + c_1 \int_0^\infty e^{- \rho t} X_t^{\xi,\eta} \circ_d d \eta_t + c_2 \int_0^\infty e^{- \rho t}X_t^{\xi,\eta} \circ_u d \xi_t \Big] < + \infty ,\\
&\E_x \Big[ \int_0^\infty e^{-\lambda t} \alpha (X_t^{\xi ,\eta} - m)^+ dt + \kappa \int_0^\infty e^{- \lambda t}  X_t^{\xi,\eta} \circ_d d \eta_t \Big] < + \infty , \\
&\P_x ( \Delta \xi_t \cdot \Delta \eta_t > 0 ) = 0 \text{ for all } t \geq 0 \text{ and } x \in \mathbb{R}_+ .
\end{split}
\end{align}
% To that end, we consider the class of controls $(\xi, \eta) \in \mathcal{U} \times \mathcal{U}$ such that \eqref{Equ: Obj 1 finite} holds.
To that end, we consider the class of controls 
$$
\mathcal{A}:=\{(\xi,\eta) \in \mathcal{U} \times \mathcal{U} :\; \eqref{Equ: Obj 1 finite} \text{ hold true} \}.
$$
For the purpose of formulating and subsequently tackling the non-zero-sum game of the government versus the legislative body, it is convenient to introduce the two sets
\begin{equation*}
\mathcal{A}_\eta := \{\xi \in \mathcal{U} :\; (\xi,\eta) \in \mathcal{A} % \mathcal{U} \times \mathcal{U} \text{ satisfies \eqref{Equ: Obj 1 finite}}
\} 
\quad \text{and} \quad 
\mathcal{A}_\xi := \{\eta \in \mathcal{U} :\; (\xi,\eta) \in \mathcal{A} %  \mathcal{U} \times \mathcal{U} \text{ satisfies \eqref{Equ: Obj 1 finite}}
\}.
\end{equation*}
% and to define $\mathcal{A}:=\mathcal{A}_\eta\times\mathcal{A}_\xi$.}
% {\color{red}[I think we need to be careful here -- it might be wrong to define $\mathcal{A}$ as above, as this creates the same issue, e.g.\ in Def.\ \ref{Def: Nash Equilibrium}, we will have ``A couple $(\xi^* , \eta^*) \in \mathcal{A}_\eta\times\mathcal{A}_\xi$", which is problematic as the Referee pointed out. 
% It might be better to do the following blue instead [What do you think? \fd{Yes, you might be right. I think it is better to proceed as you suggest (Felix)}]: \\

The problem introduced partly in \eqref{Objective: Government} and \eqref{Objective: ECB New} is therefore formulated as a non-zero-sum game between two players: 
The government (player 1) which aims at solving 
\begin{align}\label{Value Fct Government}
V_1 (x;\eta) := \inf_{\xi \in \mathcal{A}_\eta} \mathcal{J}_{x,\eta} (\xi), 
\qquad  x \in \mathbb{R}_+, 
\end{align}
for any fixed control process $\eta \in \mathcal{U}$, 
and the legislative body (player 2) which aims at solving  
\begin{align}\label{Value Fct ECB}
V_2 (x;\xi) := \inf_{\eta \in \mathcal{A}_\xi} \mathcal{I}_{x, \xi} (\eta), 
\qquad  x \in \mathbb{R}_+, 
%\hspace{1cm} \text{or} \hspace*{1cm} V_2^+ (x;\xi) := \inf_{\eta \in \mathcal{A}} \mathcal{I}_{x, \xi}^+ (\eta), \qquad x \in \mathbb{R}_+,
\end{align}
for any fixed control process $\xi \in \mathcal{U}$. 

\begin{definition}\label{Def: Nash Equilibrium}
A couple $(\xi^* , \eta^*) \in \mathcal{A}$ forms a Nash equilibrium if and only if 
\begin{align*}
\begin{cases}
    \mathcal{J}_{x,\eta^*} (\xi^*) \leq   \mathcal{J}_{x,\eta^*} (\xi)  & 
    \text{for any } \xi \in \mathcal{A}_{\eta^*}, \\
    \mathcal{I}_{x,\xi^*} (\eta^*) \, \leq   \mathcal{I}_{x,\xi^*}(\eta)  & \text{for any } \eta \in \mathcal{A}_{\xi^*}.
\end{cases}
\end{align*}
Each player's value of the game is then given by $V_1(x;\eta^*) = \mathcal{J}_{x,\eta^*}(\xi^*)$ and $V_2(x;\xi^*) = \mathcal{I}_{x,\xi^*} (\eta^*)$.
\end{definition}

The following assumptions on the model's parameters will hold true in the rest of this paper.
\begin{assumption} \label{Assumption: c1 > c2}
The model's parameters satisfy: 
\begin{itemize}
\item[(i)] $c_1 > c_2$; 
\item[(ii)] $\rho > ( p (r - g) + \tfrac{\sigma^2}{2} p (p-1))^+$, where $p$ is defined in Assumption \ref{Assumption: Cost function h};
\item[(iii)] $\lambda > r-g$;
\item[(iv)] $m > 0$.
\end{itemize}
\end{assumption}

The condition in Assumption \ref{Assumption: c1 > c2}.(i) is typically assumed in the literature on bounded-variation stochastic
control problems in order to ensure well-posedness of the optimisation problem (see, e.g., \cite{GuPh}, \cite{DeAFe} and \cite{FeRo}) and to avoid arbitrage opportunities. In economic terms, a possible interpretation is that the Keynesian multiplier is not high enough to offset the costs of liberalisation policies.

Assumption \ref{Assumption: c1 > c2}.(ii)  reflects the fact that governments are more concerned about the present than the future, since they are in power for only a limited amount of years; hence discounting future costs and benefits at a sufficiently large rate. Moreover, combining this with Assumption \eqref{Assumption: Cost function h}.(iii), the trivial policy ``never intervene on the debt ratio" is admissible, since it yields a finite expected cost. 
The latter is guaranteed also for the problem of the legislative body due to Assumption \ref{Assumption: c1 > c2}.(iii).

In this paper, we will devote our attention to the existence of Nash equilibria of the game \eqref{Value Fct Government}--\eqref{Value Fct ECB} in the class of strategies, where at least one of the players chooses a Skorokhod reflection type policy at a constant threshold. To this end, we first recall the following well known results on Skorokhod reflection.
\begin{lemma}\label{Lemma: Skorokhod: Two sided version}
Let $a,b \in \mathbb{R}_+$ with $a<b$. For  any $x \in [a,b]$ there exists a unique couple $(\xi (a), \eta (b)) \in \mathcal{A}$ that solves the Skorokhod reflection problem 
% {\color{red}[Do we have an issue with $(\xi (a), \eta (b)) \in \mathcal{A}$?]}
\begin{align*}
\tag{\textbf{SP}(a,b;x)} \text{Find } (\xi , \eta) \in \mathcal{A} \text{ s.t. } 
\begin{cases}
X_t^{\xi,\eta} \in [a,b], \P\text{-a.s.~for }t>0, \\
\int_0^T \one_{ \{X_t^{\xi,\eta} > a \}} d \xi_t = 0, \P\text{-a.s.~for any }T>0, \\
\int_0^T \one_{ \{X_t^{\xi,\eta} < b \}} d \eta_t = 0, \P\text{-a.s.~for any }T>0,
\end{cases}
\end{align*}
and it follows that supp$\{d\xi_t (a) \}\, \cap\,$supp$\{ d \eta_t (b) \} = \emptyset$. 
\end{lemma}
% In the forthcoming analysis we also use the one-sided version of this result, which we recall in the following lemma.
\begin{lemma}
For any $\eta \in \mathcal{U}$, $a \in \mathbb{R}_+$ and $x \geq a$ there exists a unique $\xi (a) \in \mathcal{A}_\eta$ solving the Skorokhod reflection problem 
\begin{align*}
\tag{\textbf{SP}($a;x,\eta$)} \text{ Find } \xi \in \mathcal{A}_\eta \text{ s.t. } 
\begin{cases}
X_t^{\xi,\eta} \geq a, \P\text{-a.s.~for }t>0, \\
\int_0^T \one_{ \{ X_t^{\xi,\eta } > a \} } d \xi_t = 0, \P\text{-a.s.~for any } T > 0. 
\end{cases}
\end{align*}
Analogously, for any $\xi \in \mathcal{U}, b \in \mathbb{R}_+$ and $x \leq b$ there exists a unique $\eta (b) \in \mathcal{A}_\xi$ solving 
\begin{align*}
\tag{\textbf{SP}($b;x,\xi$)}  \text{ Find } \eta \in \mathcal{A}_\xi \text{ s.t. } 
\begin{cases}
X_t^{\xi,\eta} \leq b, \P\text{-a.s.~for }t>0, \\
\int_0^T \one_{ \{ X_t^{\xi,\eta } < b \} } d \eta_t = 0, \P\text{-a.s.~for any } T > 0. 
\end{cases}
\end{align*}
\end{lemma}

\begin{remark}
Using the change of variable leading to \eqref{logX}, the existence of solutions to the three Skorokhod reflection problems \textbf{SP}(a,b;x), \textbf{SP}(a;x,$\eta$) and \textbf{SP}(b;x,$\xi$)  can be deduced from Proposition 2.3, Corollary 2.4 and Theorem 2.6 in \cite{Burdzy}, as well as Lemmata 2.1 and 2.2 in \cite{DeAFe2}. 
\end{remark}

Moreover, we define 
\begin{align}\label{Def: Class M}
\begin{split}
\mathcal{M} := \{ (\xi , \eta) \in \mathcal{A}: ~ &\xi \text{ solves \textbf{SP}}(a;x, \eta) \text{ or } \eta  \text{ solves \textbf{SP}}(b;x, \xi) \text{ for } x \in \mathbb{R}_+ \text{ and some } a,b \in \mathbb{R}_+ \} 
\end{split}
\end{align}
and aim to prove the existence and uniqueness of %Nash equilibria given by 
{a Nash equilibrium} $(\xi,\eta)\in \mathcal{M}$, in different parameter configurations of the game. 
Indeed, we will show that if at least one player acts according to a Skorokhod reflection type policy as specified above, the game \eqref{Value Fct Government}--\eqref{Value Fct ECB}  admits a unique Nash equilibrium.  
Clearly, when not restricting at least one of the players to a Skorokhod reflection type policy there could also exist Nash equilibria outside of the set $\mathcal{M}$. 
However, as pointed out in previous contributions such as \cite{DeAFe2}, it is impossible to rank different Nash equilibria without an additional optimality criterion.

\section{The Optimal Governmental Debt Management Rule}\label{Section: Governments problem}

In this section, we study the problem of
the government choosing their investment economic policy $\xi$, taking into account that the legislative body (e.g. Congress) may or may not choose to intervene on the debt ratio.
In the following, we distinguish between two cases, depending on the chosen control policy of the legislative body:
\begin{align}\label{Def: Control eta^b}
(\textrm{I}) ~~ \eta_t =\overline{\eta}_t := 0,  
\hspace*{2cm} (\textrm{II})~~ \eta_t = \eta^b_t := \one_{ \{t>0\}} [ (x-b)^+ + \eta_t (b) ].
\end{align}
In particular, the legislative body  does not intervene in Case (\textrm{I}), while in Case (\textrm{II}) it imposes a debt ceiling mechanism, which forces the government to keep its debt ratio below a fixed level $b \in \mathbb{R}_+$ (via e.g.~the adoption of liberalisation policies). 
In the latter definition, $\eta(b)$ uniquely solves the Skorokhod reflection problem \textbf{SP}($b; (x \wedge b), \xi$).
In the following, we aim at determining a best response (i.e.~an optimal control strategy $\xi \in \mathcal{A}_\eta$) in both cases.

For simplicity of exposition, we assume the running cost function $h(x) = x^2/2$ in \eqref{Objective: Government} (cf.\ Remark \ref{remh}) in the rest of the paper, which further yields that
\begin{equation}\label{Assumption:Rho:p=2}
\text{Assumption \ref{Assumption: c1 > c2}.(ii) with $p=2$} \quad \Leftrightarrow \quad \rho > 2(r-g)+\sigma^2.
\end{equation}

It is clearer to present the two cases (\textrm{I}) and (\textrm{II}) separately.

\subsection{The government's optimal strategy under no legislative body intervention:~Case~(\textrm{I})} \label{Sec: Government Case I}
Let us assume that the legislative body does not intervene on the government's debt. 
The value function \eqref{Value Fct Government} %\eqref{Value Fct ECB} 
thus rewrites as  
\begin{align}\label{Value function Vbar1(x)}
\overline{V}_1(x) := \inf_{\xi \in \mathcal{A}_{\overline{\eta}}} \E_x \Big[ \int_0^\infty e^{- \rho t} h(X_t^{\xi,\overline{\eta}})dt  - c_2 \int_0^\infty e^{- \rho t}X_t^{\xi,\overline{\eta}} \circ_u d \xi_t \Big]. 
\end{align}
where we let $\overline{V}_1(x):= V_1 (x;\overline{\eta})$. It follows from standard theory that we can associate the value function $\overline{V}_1$ of \eqref{Value function Vbar1(x)} with a suitable solution to the Hamilton-Jacobi-Bellman (HJB) equation
\begin{align}\label{Def: HJB Ubar1}
\min \big\{ (\mathcal{L} - \rho)u(x) + \tfrac12 x^2, ~ u'(x) - c_2 \big\} = 0
\end{align}
for all $x \in \mathbb{R}_+$, where the second order linear operator $\mathcal{L}$ defined by its action on functions $f \in C^2$ is defined by
$$
\mathcal{L} f :=\tfrac{1}{2} \sigma^2 x^2 f^{''}+(r-g)x f^{'}.
$$
We guess that the  government chooses to increase its debt ratio only when the current level is sufficiently small. 
Hence, we expect that there exists a critical $x$-level $\overline{a}$ at which the government increases their debt ratio via a Skorokhod reflection type policy. For any $x \in \mathbb{R}_+$, we thus consider the control 
\begin{align}\label{Def: Optimal Control a}
\xi_t^{\overline{a}} := \one_{\{ t>0\}} [(\overline{a} -x)^+ + \xi_t (\overline{a})]
\end{align}
where $\xi (\overline{a})$ is the unique solution to the Skorokhod reflection problem \textbf{SP}$(\overline{a}; (x \wedge \overline{a}), \overline{\eta})$. 
As a consequence, we can associate the given problem \eqref{Value function Vbar1(x)} with the free-boundary problem
\begin{align}\label{Def: Boundary problem gov case 1}
\begin{cases}
(\mathcal{L} - \rho) u(x) \geq \frac12 x^2 ,  & x \in \mathbb{R}_+, \\
(\mathcal{L} - \rho) u(x) = \frac12 x^2 , & \overline{a} < x, \\
u'(x) \geq c_2, & x \in \mathbb{R}_+, \\
u'(x) = c_2, & 0<x \leq \overline{a},\\
u''(\overline{a}) = 0, \\
\displaystyle{\lim_{x\to +\infty}} \Big(u(x)-\frac{x^2}{2(\rho-2(r-g)-\sigma^2)} \Big)=0,
\end{cases}
\end{align}
where we impose an additional smoothness condition at the boundary $\overline{a}$ and the latter one %condition in \eqref{Def: Boundary problem gov case 1} serves the purpose 
in order to guarantee uniqueness of the solution to the free-boundary problem. 
% In the following theorem, we verify that the solution to the free-boundary problem \eqref{Def: Boundary problem gov case 1} indeed coincides with the value function \eqref{Value function Vbar1(x)} and derive an optimal debt management policy for the government in Case (I). The proof can be found in Appendix \ref{Appendix:ProofVTG1}.
In the following theorem, we derive a solution to the optimal debt management problem of the government in Case (\textrm{I}). 
We follow the usual \textit{guess-and-verify} approach by first constructing a solution to the free-boundary problem  \eqref{Def: Boundary problem gov case 1} and then verify that this \textit{candidate} value function indeed coincides with the true value function \eqref{Value function Vbar1(x)}. As a byproduct, we obtain the optimal debt management policy, which prescribes to reflect the debt ratio upwards according to the Skorokhod reflection type policy \eqref{Def: Optimal Control a}. The proof can be found in Appendix \ref{Appendix:ProofVTG1}.

\begin{theorem}[Verification Theorem: Case (\textrm{I})]\label{VT:Gcase1}
Assume that the legislative body does not intervene on the debt ratio and thus acts according to the policy $\overline{\eta} \equiv 0$. 
Then, the value function $\overline{V}_1$ of \eqref{Value function Vbar1(x)} is given by
\begin{align}
\overline{V}_1 (x) := 
\begin{cases}
\overline{V}_1 (\overline{a}) - c_2 (\overline{a} -x), & 0 < x \leq \overline{a} , \\
\overline{D}_1 (\overline{a}) x^{\delta_2} + \frac{1}{2(\rho - 2(r-g)-\sigma^2)} x^2 , & \overline{a} < x ,
\end{cases}
\end{align} 
where 
\begin{align}\label{Def: Dbar1 and a*1}
\overline{D}_1 (a) := - \frac{1}{(\rho - 2(r-g)-\sigma^2)\delta_2 (\delta_2 -1) a^{\delta_2 -2}}, ~~~~\text{and} ~~~~ \overline{a} := \frac{(1-\delta_2 )c_2(\rho-2(r-g)-\sigma^2)}{(2-\delta_2 )},
\end{align}
with $\delta_2$ denoting the negative root to the equation $\frac{1}{2} \sigma^2 \delta (\delta-1) + (r-g)\delta - \rho = 0$. Moreover, the
admissible $\xi_t^{\overline{a}}$ of \eqref{Def: Optimal Control a}, with $\overline{a}$ given by \eqref{Def: Dbar1 and a*1}, is optimal for problem \eqref{Value function Vbar1(x)}.
\end{theorem}

\subsection{The government's optimal strategy under legislative body interventions:~Case~(\textrm{II})} \label{Sec: Government Case II}
We begin by fixing a constant $b \in \mathbb{R}_+$ and assume that the legislative body acts according to the control policy $\eta^b$ of \eqref{Def: Control eta^b}, thus keeping the debt ratio below the debt ceiling $b$ according to a Skorokhod reflection type policy. From the government's point of view, we thus study the problem $V_1(x;\eta^b)$ defined in \eqref{Value Fct Government} and given by\footnote{\label{foot3}For ease of notation, we denote by $V_1(x;y)$ and $V_2(x;y)$, $x,y \in \mathbb{R}_+$, the control value functions $V_1(x;\eta^y)$ and $V_2(x;\xi^y)$, i.e. when their opponents choose the Skorokhod reflection type strategies $\eta^y$ and $\xi^y$, respectively.}
\begin{align}\label{Value function V1(x,b)}
V_1(x;b) := \inf_{\xi \in \mathcal{A}_{\eta^b}} \E_x \Big[ \int_0^\infty e^{- \rho t} h(X_t^{\xi,\eta^b})dt + c_1 \int_0^\infty e^{- \rho t} X_t^{\xi, \eta^b} \circ_d d \eta_t^b - c_2 \int_0^\infty e^{- \rho t}X_t^{\xi,\eta^b} \circ_u d \xi_t \Big]. 
\end{align}
%\nr{where we let $V_1(x;b) := V_1 (x;\eta^b)$.[I don't like this mix of notations - created a footnote, check]} 
Again, we can associate the latter value function with the solution to the HJB equation 
\begin{align}\label{Def: HJB U1(x,b)}
\min \big\{ (\mathcal{L}-\rho) u(x;b) + \tfrac12 x^2, ~ u' (x;b) -c_2 \big\} = 0
\end{align}
for all $x \in (0,b)$ with boundary condition $u(0;b) = 0$ and Neumann boundary condition $u'(b;b) = c_1$. 

We guess that the government increases their debt ratio only when the current level is sufficiently small. 
Hence, we expect that for any given debt ceiling $b \in \mathbb{R}_+$, there exists a critical debt-issuance level $a(b)$ at which the government increases their debt ratio with \textit{minimal effort}, via a Skorokhod reflection type policy, where we stress the (possible) dependency on the debt ceiling threshold $b \in \mathbb{R}_+$. For any $x \in \mathbb{R}_+$, we thus consider the control 
\begin{align}\label{Def: Optimal Control a(b)}
\xi^{a(b)}_t := \one_{ \{ t > 0 \}} [ (a(b) -x)^+ + \xi_t (a(b)) ],
\end{align} 
where $\xi (a(b))$ 
is the unique control such that the couple $(\xi(a(b)),\eta(b))$ 
solves the Skorokhod reflection problem \textbf{SP}($a(b),b; (x \vee a(b)) \wedge b$). 
As a consequence, we can associate the given problem \eqref{Value function V1(x,b)} %\eqref{Value Fct Government} 
with the free-boundary problem
\begin{align}\label{Free Boundary Problem Government}
\begin{cases}
(\mathcal{L}- \rho)u(x;b) \geq - \frac12 x^2, & 0< x <b, \\
(\mathcal{L}- \rho)u(x;b) = - \frac12 x^2, & a(b) < x< b, \\
u'(x;b) \geq c_2, &0<x < b, \\
u'(x;b) = c_2, & 0< x \leq a(b), \\
u'(x;b) = c_1, & b \leq x, \\
u''(a(b);b) = 0,
\end{cases}
\end{align}
where we impose an additional smoothness condition at the free boundary $a(b)$. 
The forthcoming analysis is dedicated to determining the optimal debt-issuance threshold $a(b)$ and proving the optimality of the control \eqref{Def: Optimal Control a(b)} for the original debt ratio management problem of the government \eqref{Value function V1(x,b)}, which corresponds to \eqref{Value Fct Government} with $\eta = \eta^b$ defined in \eqref{Def: Control eta^b}.

We begin with solving the free-boundary problem \eqref{Free Boundary Problem Government} by %following a \textit{guess-and-verify approach}. By 
constructing a solution to the ordinary differential equation and imposing the boundary conditions to obtain a \textit{candidate} value function 
\begin{align}\label{Def: Cand VF V1}
U_1(x;b) = \begin{cases}
U_1(a(b);b) - c_2 (a (b) - x), & 0 < x \leq a(b), \\
D_1(a(b)) x^{\delta_1} + D_2(a(b)) x^{\delta_2} + \frac{1}{2(\rho -2 (r-g) -\sigma^2)} x^2, & a(b) < x < b, \\
U_1(b;b) + c_1 (x-b) , & b \leq x,
\end{cases}
\end{align}
with 
\begin{align*}
D_i (a) = \frac{(\delta_{3-i} -2)a - c_2(\delta_{3-i} -1)(\rho -2(r-g)-\sigma^2)}{(-1)^{i+1} \delta_i(\delta_1 - \delta_2)(\rho - 2(r-g) -\sigma^2)a^{\delta_i -1}}, \qquad \text{for } i=1,2,
\end{align*}
and the constants $\delta_1, \delta_2$ denoting the positive and negative roots to the equation $\frac{1}{2} \sigma^2 \delta (\delta-1) + (r-g)\delta - \rho = 0$, respectively. 
The optimal boundary is given by the solution $ a(b) \in (0,b)$ to the equation
\begin{align}\label{Equ: F = 0}
F(a(b),b) = 0,
\end{align}
where we define
\begin{align}\label{Def: F(a,b)}
F(a,b) := &[(2- \delta_2)a - c_2 (1 - \delta_2)(\rho - 2 (r-g)-\sigma^2)]\Big(\frac{b}{a}\Big)^{\delta_1 -1} \nonumber\\
+&[(\delta_1 -2)a -c_2(\delta_1 -1)(\rho - 2(r-g)-\sigma^2)]\Big(\frac{b}{a} \Big)^{\delta_2 -1} %\nonumber \\&
-(\delta_1 -\delta_2)[b -c_1(\rho -2(r-g)-\sigma^2)]. %\nonumber
\end{align}

In the following result, we prove the existence and uniqueness of the optimal boundary, as well as some monotonicity and limit properties that will be useful in the search for a Nash equilibrium later in Section \ref{nash}. Its proof can be found in Appendix \ref{Appendix:ProofLemmaBounda(b)}. 

\begin{lemma}\label{Lemma: a(b) unique}
Let $b \in \mathbb{R}_+$ and recall $F(\cdot,b)$ defined by \eqref{Def: F(a,b)} on $(0,b)$. 
Then, \\[1mm]
{\rm (i)} There exists a unique $a(b) \in (0, b)$ solving %$a(b) \in (0, \tilde{a} \wedge b)$ that satisfies
$F(a(b),b) = 0$ in \eqref{Equ: F = 0}, that satisfies $\tfrac{\partial}{\partial a} F(a(b),b) > 0$ and $a(b) < \tilde{a}$, where 
\begin{align}\label{Def: atilde}
 \tilde{a}:= 
% \frac{c_2(\delta_1 -1)(1-\delta_2)(\rho-2(r-g)-\sigma^2)}{(\delta_1-2)(2-\delta_2)} =
c_2 (\rho - (r-g)) > 0.
\end{align}
{\rm (ii)} We have 
\begin{equation}\label{a'0-inf}
\text{$a(\cdot)$ is} 
\begin{cases} 
\text{increasing on $(0, \hat b)$}, \\
\text{decreasing on $(\hat b, \infty)$}, 
\end{cases}
\quad \text{where $\hat{b}$ is the unique solution to $\frac{\partial}{\partial b} F (a(\hat{b}),\hat{b})=0$},
\end{equation}
and $\lim_{b\to 0} a(b) = 0$ as well as $\lim_{b\to \infty} a(b) = \overline{a}$, where $\overline{a}$ is the optimal debt-issuance threshold defined by \eqref{Def: Dbar1 and a*1} in Case (I) of non-intervention by a legislative body. \\[1mm]
{\rm (iii)} Furthermore, $b \mapsto a(b)$ is concave on the interval $(0,\hat{b})$. 
\end{lemma}

% The next result presents the solution to problem \eqref{Value function V1(x,b)} and the optimality of the control $\xi^{a(b)}$. 
In the following theorem, we prove that the value function $V_1$ of \eqref{Value function V1(x,b)} is indeed given by the solution to the free-boundary problem \eqref{Free Boundary Problem Government}, i.e.~the function $U_1$ of \eqref{Def: Cand VF V1}.  Moreover, the control policy $\xi^{a(b)}$ of \eqref{Def: Optimal Control a(b)}, with $a(b)$ given as the unique solution to equation \eqref{Equ: F = 0}, is proven to be optimal.
Before doing so, we first notice that the control policy $\xi^{a(b)}$ as in \eqref{Def: Optimal Control a(b)}, combined with the policy $\eta^b$ of \eqref{Def: Control eta^b}, is indeed admissible. 
Clearly, the couple solves \textbf{SP}($a(b),b;x$) and as such belongs to $\mathcal{A}$. 
Indeed, by arguing as in Lemma 4.1 in \cite{ShLeGa} one can easily show \eqref{Equ: Obj 1 finite} and moreover, $\P_x (\Delta \xi^{a(b)} \cdot \Delta \eta^b > 0) = 0$ for all $t \geq 0$, by construction. The proof of the following result, which concludes this section, can be found in Appendix \ref{Appendix:ProofVTG2}.

\begin{theorem}[Verification Theorem: Case (\textrm{II})]\label{VT:Gcase2}
Assume that the legislative body acts according to the control policy $\eta^{b}$ of \eqref{Def: Control eta^b}. Then, the function $U_1$ of \eqref{Def: Cand VF V1} coincides with the government's value function $V_1$ in \eqref{Value function V1(x,b)} %{Value Fct Government} 
and the admissible $\xi^{a(b)}$ of \eqref{Def: Optimal Control a(b)} is optimal for problem \eqref{Value function V1(x,b)}. %{Value Fct Government}.
\end{theorem}

\section{The optimal debt ceiling}\label{Section: The institutions problem}
In this section, we study the control problem of the legislative body. 
As seen in Section \ref{Section: Governments problem}, the best response of the government to either a legislative body non-intervention policy, or a debt ceiling mechanism (threshold-type policy) %imposed by the legislative body 
is given by a debt-issuance threshold-type policy. 
We now reverse the roles and assume that the government chooses to increase its debt ratio at a certain level $a \in \mathbb{R}_+$, i.e.~to the debt-issuance control policy 
\begin{align}\label{Def: Control xi a}
\xi^a_t := \one_{\{ t > 0 \} } [ (a-x)^+ + \xi_t (a) ],
\end{align} 
where $\xi (a)$ uniquely solves the Skorokhod reflection problem \textbf{SP}($a; x \vee a, \eta$). 
In the following, for any such level $a$, we study the problem \eqref{Value Fct ECB} of finding a best response (i.e.~an optimal control strategy $\eta \in \mathcal{A}_{\xi^a}$). 
%To keep notation short, we let $V_2(x;a):= V_2(x;\xi^a)$,
%as well as $V_2^+(x;a):= V_2^+(x;\xi^a)$ and 
We thus consider the problem$^{\ref{foot3}}$ 
\begin{align}\label{Value function V2 (x,a)}
V_2(x;a) &:= \inf_{\eta \in \mathcal{A}_{\xi^a}} \E_x \Big[ \int_0^\infty e^{- \lambda t} \alpha(X_t^{\xi^a,\eta} - m)^+ dt +   \kappa \int_0^\infty e^{-\lambda t} X_t^{\xi^a,\eta} \circ_d d \eta_t \Big]. 
%\\ V_2^+(x;a) &:= \inf_{\eta } \E_x \Big[ \int_0^\infty e^{- \lambda t} \alpha(X_t^{\xi^a,\eta} - m)^+ dt +   \kappa \int_0^\infty e^{-\lambda t} X_t^{\xi^a,\eta} \circ_d d \eta_t \Big].\label{Value ECB New}
\end{align}
Via standard arguments, we can associate the value function $V_2$ of \eqref{Value function V2 (x,a)} with a suitable solution to the HJB equation 
\begin{align}\label{Def: HJB ECB}
\min \big\{ (\mathcal{L} - \lambda ) u(x;a) + \alpha (x-m)^+ , \kappa - u'(x;a) \big\} &= 0, 
\end{align}
for all $x \in (a,\infty)$ with Neumann boundary condition $u'(a;a) = 0$. 
%In the forthcoming, we follow a \textit{guess-and-verify approach} and 
We presume that the legislative body may only decrease the debt ratio when the current level is sufficiently large. 
Therefore, if the legislative body chooses to intervene, we expect that for any given debt-issuance threshold $a \in \mathbb{R}_+$, there exists a critical debt ceiling level $b(a)$ at which the legislative body forces a decrease in the debt ratio via a Skorokhod reflection type policy, where we stress the (possible) dependency on the debt-issuance threshold $a \in \mathbb{R}_+$. 
On the other hand, also a non-intervention policy is conceivable. 
As it turns out, it is crucial in our analysis to distinguish two different cases, depending on the legislative body's time preference rate $\lambda$:
\begin{align*}
(\textrm{I}) ~~ \lambda > r-g + \frac{\alpha}{\kappa},  \vspace*{0.2cm}
\hspace*{2cm} (\textrm{II})~~ \lambda < r-g + \frac{\alpha}{\kappa}.
\end{align*}
In the forthcoming Sections \ref{Sec: ECB Case I} and \ref{Sec: ECB Case II} we study these cases separately, providing an optimal control strategy by the legislative body for each one of them. 

\subsection{The legislative body's optimal strategy under high time preference rate:~Case~(\textrm{I})}
\label{Sec: ECB Case I}
Notice that the legislative body discounts future events with a relatively large discount factor in this case, and it is therefore appropriate to assume that the legislative body disregards the risk of future government insolvency at a greater extent compared to Case (\textrm{II}). 
We verify this intuition by showing that indeed, the optimal control policy of the legislative body prescribes \textit{not to intervene} on the debt ratio at all. 
%This result holds true for both running cost functions. 

To this end, we prove that the value function $V_2$ of \eqref{Value function V2 (x,a)} coincides with a suitable solution to a fixed-boundary problem 
\begin{align}\label{Def: Free Boundary Problem ECB b infty}
\begin{cases}
(\mathcal{L}-\lambda ) u (x;a) = - \alpha (x-m)^+, & a < x, \\
u' (x;a) < \kappa, & a < x, \\
u'(x;a) = 0, & 0 < x \leq a.
\end{cases}
\end{align}
We can solve the latter problem by  constructing a solution to the ordinary differential equation and imposing the stated boundary condition. 
% In the following theorem, we verify that the solution to the fixed-boundary problem \eqref{Def: Free Boundary Problem ECB b infty} indeed coincides with the value function $V_2$ of \eqref{Value function V2 (x,a)}. Its proof can be found in Appendix \ref{Appendix:ProofVTI1}.
In the following theorem, we again follow the usual \textit{guess-and-verify} approach by proving that the constructed solution to the fixed-boundary problem \eqref{Def: Free Boundary Problem ECB b infty} indeed coincides with the value function $V_2$ of \eqref{Value function V2 (x,a)}. Moreover, the proof reveals that the no-intervention policy $\overline{\eta} = 0$ is optimal. The proof, including the construction of the solution to the fixed-boundary problem \eqref{Def: Free Boundary Problem ECB b infty}, can be found in Appendix \ref{Appendix:ProofVTI1}.

\begin{theorem}[Verification Theorem: Case (\textrm{I})] \label{VT:Icase1}
Assume that the government acts according to the control policy $\xi^a$ of \eqref{Def: Control xi a}. 
Then, the function $V_2$ of \eqref{Value function V2 (x,a)} is given by 
\begin{align}\label{Def: Candidate U2 b infty}
%\overline{U}_2 (x;a) &= \begin{cases}\overline{U}_2 (a;a) & 0 < x \leq a, \\ \overline{D}_2 (a) x^{\theta_2} + \alpha (\frac{x}{\lambda - (r-g)} - \frac{m}{\lambda} ), & a< x,\end{cases} \\
%\label{Def U2 b infty new}
%\overline{U}_2 (x;a) 
V_2 (x;a) &= 
\begin{cases} 
V_2 (a;a), & 0 < x \leq a; \\
\overline{D}_2 (a) x^{\theta_2} + H(x),\hspace*{1.7cm} & a<x, 
\end{cases}
\end{align}
where %\nr{[Is my bracket below correct? Also is $d(x)$ the same as $d(x,t)$?]}\fd{[Yes and yes. I will stick to the notation $d(x,t)$, to show the time dependence]}
\begin{align}
%\overline{D}_2 (a) := - \frac{\alpha}{\theta_2 (\lambda - (r-g))} a^{1-\theta_2}, ~~~~~~~~
\overline{D}_2(a) &:= - \frac{\alpha}{ \theta_2} a^{1 -\theta_2 } \int_0^\infty e^{-(\lambda - (r-g))t} \Phi (d_1(a,t))dt , \nonumber \\
H(x) &:= \alpha \int_0^\infty \Big ( xe^{-(\lambda -(r-g))t} \Phi (d_1(x,t)) - m e^{ - \lambda t} \Phi(d_2(x,t)) \Big )dt, \label{Def:FunctionsforVTI1}\\
d_1 (x,t) &:= \frac{\log \big(\frac{x}{m}\big) + (r-g+\frac12 \sigma^2)t}{\sigma \sqrt{t}} \hspace{0.5cm} \text{and} \hspace*{0.5cm}   d_2(x,t) := d_1(x,t) - \sigma \sqrt{t}, \nonumber
\end{align}
with $\theta_2$ denoting the negative root to the equation $\frac12 \sigma^2 \theta (\theta - 1) + (r-g) \theta - \lambda = 0$ and $\Phi(\cdot)$ denoting the cumulative distribution function of a standard normal random variable. The optimal policy for the legislative body prescribes not to act on the debt ratio, i.e. $\overline{\eta} := 0$. 
\end{theorem}

%\begin{remark}\label{Rem: a infty} 
%\fd{[Have another look: We can leave it here (to already anticipate the equilibrium in this case) or postpone it to section 5]}
Notice that the strategy of the legislative body not to intervene on the debt ratio is not triggered by some specific $a \in \mathbb{R}_+$. 
This comes solely from the fact that $\lambda > r-g+\frac{\alpha}{\kappa}$, independently of the governmental choice of a debt-issuance level $a \in \mathbb{R}_+$. 
On the other hand, the best response of the government to a legislative body non-intervention policy $\overline{\eta}$ from \eqref{Def: Control eta^b}, has been treated in Section \ref{Sec: Government Case I}. 
In particular, the optimal control for problem $V_1(x; \overline{\eta})$ of \eqref{Value Fct Government} (cf.\ $\overline{V}_1(x)$ in \eqref{Value function Vbar1(x)}) is given by the Skorokhod reflection type policy $\xi^{\overline{a}}$ as in \eqref{Def: Optimal Control a(b)}, with $\overline{a}$ defined in \eqref{Def: Dbar1 and a*1}.
Indeed, we prove in Section \ref{nash} that this pair of strategies leads to an equilibrium.

\subsection{The legislative body's optimal strategy under low time preference rate:~Case~(\textrm{II})}\label{Sec: ECB Case II}
While it is optimal for the legislative body to never intervene in Case (\textrm{I}), we show that in this case the best response to a governmental policy \eqref{Def: Control xi a} requires intervention. 
Indeed, for any given governmental debt-issuance threshold $a \in \mathbb{R}_+$, this will prescribe keeping the debt ratio below a certain debt ceiling $b(a)$ with \textit{minimal effort}, via a Skorokhod reflection type policy, where we stress the (possible) dependency on the debt-issuance threshold $a \in \mathbb{R}_+$. 
%%%
For any $x \in \mathbb{R}_+$, we thus consider the control
\begin{align}\label{Def: Control eta b(a)}
\eta_t^{b(a)} := \one_{ \{ t > 0\} } [(x - b(a))^+ + \eta_t (b(a)) ] ,
\end{align}
where $\eta(b(a))$ is the unique control such that the couple $(\xi (a), \eta (b(a)))$ solves the Skorokhod reflection problem \textbf{SP}($a,b(a); (x \vee a) \wedge b(a)$). 
%%%
As a consequence, we can associate the given problem \eqref{Value function V2 (x,a)} with the free-boundary problem
\begin{align}\label{Def: Free Boundary problem ECB}
\begin{cases}
(\mathcal{L} -\lambda) u(x;a) \geq - \alpha (x-m)^+, & a < x, \\
(\mathcal{L} -\lambda) u(x;a) = - \alpha (x-m)^+, & a < x < b(a), \\
u'(x;a) \leq \kappa, & a < x, \\ %\text{\nr{[Should this be for $x < b(a)$ instead?]}}
u'(x;a) = \kappa, &b(a) \leq x, \\
u'(x;a) = 0, & 0 < x \leq a, \\
u''(b(a);a) = 0,
\end{cases}
\end{align}
where we imposed an additional smoothness condition at the free boundary $b(a)$. 
The forthcoming analysis is dedicated to determining the optimal threshold $b(a)$ and proving the optimality of the control \eqref{Def: Control eta b(a)} for the original debt ratio management problem of the government \eqref{Value function V2 (x,a)}, which corresponds to \eqref{Value Fct ECB} with $\xi = \xi^a$ defined in \eqref{Def: Control xi a}.

We begin with solving the free-boundary problem \eqref{Def: Free Boundary problem ECB} by %following a \textit{guess-and-verify approach}. By 
constructing a solution to the ordinary differential equation and imposing the boundary conditions to obtain a \textit{candidate} value function

\begin{align}\label{Def: Candidate VF U2}
%U_2 (x;a) = \begin{cases} U_2 (a;a), & 0 < x < a, \\ D_3(b(a)) x^{\theta_1} + D_4(b(a)) x^{\theta_2} + \alpha \big( \frac{x}{\lambda -(r-g)} - \frac{m}{\lambda}\big), & a < x < b(a) \\ U_2(b(a);a) + \kappa (x - b(a)), &  b(a) < x, \end{cases} \end{align} as well as  \begin{align}\label{Def U2+}
U_2 (x;a) := 
\begin{cases}
U_2 (a;a), & 0<x\leq a, \\
D_3 (b(a)) x^{\theta_1} + D_4 (b(a)) x^{\theta_2} + H(x),\hspace*{1cm} &a<x<b(a), \\
U_2 (b(a));a) + \kappa (x - b (a)), & b(a) \leq x,
\end{cases}
\end{align}
where 
\begin{align*}
%D_i (b) &:= \frac{(k(\lambda -(r-g)) -\alpha)(\theta_{5-i} -1)}{(-1)^{i} \theta_{i-2} (\theta_1- \theta_2)(\lambda - (r-g)) b^{\theta_{i-2} -1}}, \qquad \\
D_i(b) 
&:= \frac{ b^{ 1- \theta_{i-2}}}{\theta_{i-2} (\theta_2 - \theta_1)}\Big[(\theta_{5-i} -1)\Big( k -\alpha \int_0^\infty e^{-(\lambda - (r-g))t} \Phi (d_1(b,t)) dt \Big)\\
&\hspace*{6cm} + \alpha  \int_0^\infty e^{-(\lambda -(r-g))t} \frac{1}{\sqrt{2 \pi t}\sigma} e^{- \frac{1}{2} d_1(b,t)^2} dt \Big], \qquad  i = 3,4,
\end{align*}
and the constants $\theta_{1}, \theta_{2}$ are given by the positive and negative roots to the equation $\frac12 \sigma^2 \theta (\theta - 1) + (r-g) \theta - \lambda = 0$, respectively.
The optimal boundary is given by the solution $b(a) \in (a,\infty)$ to the equation 
\begin{align}\label{Equ: G=0}
G(a,b(a)) = 0, 
%\hspace*{1.5cm} G^+ (a, b^+(a)) = 0
\end{align}
where $G(a,\cdot)$ is defined on $(a,\infty)$ by 
\begin{align} \label{Def: G(a,b)}
G (a,b) &= \Big[(\theta_1-1)\Big( \frac{b}{a}\Big)^{1- \theta_2} + (1-\theta_2)\Big( \frac{b}{a}\Big)^{1 - \theta_1} \Big] \Big( \frac{\kappa}{\theta_1 - \theta_2} - \frac{\alpha}{(\theta_1- \theta_2)(\lambda -(r-g))} \one_{\{b \geq m\}} \Big) \nonumber\\ 
&~~~~+\frac{\alpha}{(\theta_1-\theta_2)(\lambda - (r-g))} \Big[ (1 - \theta_2)  \Big( \frac{m}{a}\Big)^{1 - \theta_1}  \one_{\{ b>m>a\}}  + (\theta_1 -1) \Big( \frac{m}{a}\Big)^{1 - \theta_2}  \one_{\{b>m > a\}} \Big] \nonumber \\
&~~~~+\frac{\alpha}{(\lambda - (r-g)} \one_{\{a\geq m\}}.
% &= \frac{\alpha}{(\theta_1- \theta_2)(\lambda -(r-g))} \Big[ (\theta_1 -1)  \Big( \Big( \frac{m}{a}\Big)^{1- \theta_2} - \Big( \frac{b}{a}\Big)^{1 - \theta_2} \Big) - (1- \theta_2) \Big( \frac{b}{a}\Big)^{1- \theta_1} \Big] \one_{\{b > m\}} \nonumber\\
% &~~~~+ \frac{\alpha}{(\theta_1-\theta_2)(\lambda - (r-g))} \Big[ (1 - \theta_2) \nr{\Big( [?]}
% \Big( \frac{m}{a}\Big)^{1 - \theta_1} \nr{\one_{\{ a<m \leq b\}}} %[ \one_{\{ m>a\}} - \one_{\{ m>b\}} ] 
% + (1 - \theta_1) \Big( \frac{m}{a}\Big)^{1 - \theta_2} \one_{\{a>m\}}\Big] \nonumber\\
% &~~~~+ \Big[(\theta_1-1)\Big( \frac{b}{a}\Big)^{1- \theta_2} + (1-\theta_2)\Big( \frac{b}{a}\Big)^{1 - \theta_1} \Big] \frac{\kappa}{\theta_1 - \theta_2} + \frac{\alpha}{(\lambda - (r-g)} \one_{\{a\geq m\}}. 
\end{align}

In the following lemma we state our results on the existence and uniqueness of a solution $b(a) \in (a, \infty)$ solving \eqref{Equ: G=0}, as well as some monotonicity and limit properties that will be useful in the search for a Nash equilibrium later in Section \ref{nash}. Its proof can be found in Appendix \ref{Appendix:ProofLemmaBoundb(a)}.

\begin{lemma} \label{Lem:b(a)}
Let $a \in \mathbb{R}_+$ and recall $G(a,\cdot)$ defined by \eqref{Def: G(a,b)} on $(a,\infty)$. Then, \\[1mm]
{\rm (i)} There exists a unique $b(a) \in (a, \infty)$ %$b(a) \in (a \vee m, \infty)$ 
solving $G(a, b(a)) = 0$ in \eqref{Equ: G=0}, that satisfies $\frac{\partial}{\partial b} G(a, b(a))<0$. \\[1mm]
{\rm (ii)} We have  
\begin{align}\label{Def: b0}
 b (a) \geq b_0 := \Big( \frac{\alpha}{\alpha - \kappa(\lambda-(r-g))}\Big)^{\frac{1}{1-\theta_2}}m > m.
\end{align}
{\rm (iii)} The function $a \mapsto b(a)$ is strictly increasing on $\mathbb{R}_+$. 
In particular, it takes a linear form $b(a) = (1/\tilde{q}) a$, for all $a >m$, where $\tilde{q} \in (0,1)$ is given by the solution to
\begin{align}\label{Def: tilde q}
    (1 - \theta_2) ( \kappa (\lambda - (r-g)) - \alpha ) \tilde{q}^{\,\theta_1 - 1} + (\theta_1 -1) (\kappa (\lambda - (r-g)) - \alpha) \tilde{q}^{\,\theta_2 -1} + \alpha (\theta_1 - \theta_2) = 0. 
\end{align}
{\rm (iv)} Moreover, $a \mapsto b(a)$ is convex on the interval $(0,m)$, with $\lim_{a\to 0} b(a) = b_0$, where $b_0$ is given by \eqref{Def: b0}, and $\lim_{a\to \infty} b(a) = \infty$.
\end{lemma}

Before we present the optimality of the controls, we first notice that the control policy $\eta^{b(a)}$ as in \eqref{Def: Control eta b(a)}, combined with the policy $\xi^a$ of \eqref{Def: Control xi a}, is indeed admissible. 
Clearly, the couple solves \textbf{SP}($a,b(a);x$) and as such belongs to $\mathcal{A}$. 
Indeed, by arguing as in Lemma 4.1 in \cite{ShLeGa} one can easily show \eqref{Equ: Obj 1 finite} and moreover, $\P_x (\Delta \xi^{a} \cdot \Delta \eta^{b(a)} > 0) = 0$ for all $t \geq 0$, by construction.
% The proof of the next result can be found in Appendix \ref{Appendix:ProofVTI2}.
In the next theorem, we conclude this chapter by proving that the value function \eqref{Value function V2 (x,a)} is indeed given by the solution to the free-boundary problem \eqref{Def: Free Boundary problem ECB}. Moreover, it is revealed that the control policy $\eta^{b(a)}$ of \eqref{Def: Control eta b(a)} is optimal for the legislative body in Case (\textrm{II}).

\begin{theorem}[Verification Theorem: Case (\textrm{II})] \label{VT:Icase2}
Assume that the government acts according to the control policy $\xi^a$ of \eqref{Def: Control xi a}. Then, the function $U_2$ of \eqref{Def: Candidate VF U2} coincides with the value function $V_2$ of \eqref{Value function V2 (x,a)}. Furthermore, the policy $\eta^{b(a)}$ of \eqref{Def: Control eta b(a)} with the optimal threshold determined via \eqref{Equ: G=0} is optimal. 
\end{theorem}

\section{Nash Equilibria in the model}\label{nash}

Our main results concern the existence of Nash equilibria in our model. These results stem from the analysis of the decision problems faced by the government and the legislative body developed in the previous Sections \ref{Section: Governments problem} and \ref{Section: The institutions problem}, respectively. 
We will focus on the two cases-- each player's best response to a Skorokhod-reflection type strategy is either a Skorokhod-reflection type strategy or a no-intervention policy -- which suggest that we should aim at determining a Nash equilibrium via its Definition \ref{Def: Nash Equilibrium} in the class $\mathcal{M}$ of \eqref{Def: Class M}, in which at least one player acts according to a Skorokhod-reflection type policy. The analysis in this section focuses on this direction.

We highlight the peculiarity arising from our results in Section \ref{Section: The institutions problem}, where we show that the legislative body may choose {\it not to intervene} on the debt ratio at all. 
Interestingly, we prove that the optimality of adopting such a strategy relies solely on their (individual) time preferences compared to the parameter constellation in the model -- it does {\it not depend} on the actions of the opposing player (government) -- see specifically our results in Section \ref{Sec: ECB Case I}. 
We thus split our search for Nash equilibria in the forthcoming analysis based on the magnitude of the legislative body time preferences.

\subsection{The case of $\lambda > r-g+\alpha/\kappa$} 
\label{sec:largelambda}

Our results in Section \ref{Sec: ECB Case I}, suggest that the legislative body should restrain themselves 
from reflecting the debt ratio at any threshold, when their time preference rate $\lambda$ is relatively large. 
In light of this non-intervention policy, it is natural to then examine what should the governmental strategy be as a best response. 
The characterisation of such a strategy is in fact the main aim of Section \ref{Sec: Government Case I}, which studies the optimal control (debt issuance policy) of the government when they are the sole player (there is no opponent).  

We present the resulting Nash equilibrium in the following theorem. 
Its proof is a simpler version of the one for Theorem \ref{Theorem: Existence and Uniqueness of NE}, and it is thus omitted for brevity.

\begin{theorem}[Existence and Uniqueness of Nash Equilibrium: Case (I)] \label{Theorem:NE1}
Suppose that the model's parameters satisfy Assumptions \ref{Assumption: c1 > c2} as well as $\lambda > r-g+\alpha/\kappa$. 
A unique Nash equilibrium of the game \eqref{Value Fct Government}--\eqref{Value Fct ECB} in the set $\mathcal{M}$ of \eqref{Def: Class M} can be characterised by the couple of controls $(\xi^{\overline{a}}, \overline{\eta}) = (\xi^{\overline{a}}, 0)$, with the former component defined as in \eqref{Def: Optimal Control a} and the threshold $\overline{a}$ is given explicitly by \eqref{Def: Dbar1 and a*1}.
\end{theorem}

\subsection{The case of $\lambda \in (r-g, r-g + \alpha/\kappa)$} 
\label{sec:smalllambda}
Our results from Sections \ref{Section: Governments problem} and \ref{Section: The institutions problem} on each player's best response 
suggest that a Nash equilibrium could be characterised by Skorokhod-reflection type policies at finite thresholds. 
More precisely, while the government increases the debt ratio at $a(b)$ (as a best response to a debt ceiling $b \in \mathbb{R}_+$), the legislative body forces a debt ratio reduction at a debt ceiling $b(a)$ (as a best response to a governmental debt-issuance threshold $a \in \mathbb{R}_+)$. 

The aim of the following theorem is to first prove that there always exists a pair $(a^*,b^*)$ forming a fixed point of these best-response-maps, such that $a^*=a(b^*)$ and $b^*=b(a^*)$, and second, that this pair is unique in the set $\mathcal{M}$, in which at least one player plays a Skorokhod-reflection type policy. 
Its proof is given in Appendix \ref{Appendix:NashEx}.

\begin{theorem}[Existence and Uniqueness of Nash Equilibrium: Case (II)] 
\label{Theorem: Existence and Uniqueness of NE}
Suppose that the model's parameters satisfy Assumptions \ref{Assumption: c1 > c2} as well as $\lambda \in (r-g, r-g+\alpha/\kappa)$. 
A Nash equilibrium in the game \eqref{Value Fct Government}--\eqref{Value Fct ECB} can be characterised by the couple of Skorokhod-reflection type policies $(\xi^{a^*}, \eta^{b^*})$, as defined in \eqref{Def: Optimal Control a(b)} and \eqref{Def: Control eta b(a)}, respectively. 
The pair of thresholds $(a^*,b^*) \in \mathbb{R}_+^2$ solves the coupled system of equations $F(a^*,b^*) = 0 = G(a^*,b^*)$ and satisfies $a^* < b^*$,  $a^* = a(b^*) < b^* = b(a^*)$ according to Lemmata \ref{Lemma: a(b) unique} and \ref{Lem:b(a)},  for a unique $b^*>b_0$ where the latter is defined in \eqref{Def: b0}.
Moreover, the Nash Equilibrium is unique in the class $\mathcal{M}$ specified in \eqref{Def: Class M}.
\end{theorem}

The following corollary reveals the connection between the equilibria determined in Theorem \ref{Theorem:NE1}-\ref{Theorem: Existence and Uniqueness of NE} by considering the limit as $\lambda \uparrow r-g+\alpha / \kappa.$ Its proof can be found in Appendix \ref{Appendix:Coro}.

\begin{corollary}\label{Coro:ConvNash}
    Let $(a^* (\lambda), b^* (\lambda)) = (a^*,b^*)$ denote the pair of thresholds that characterises the unique Nash equilibrium derived in Theorem \ref{Theorem: Existence and Uniqueness of NE}, where we stress the dependency on the parameter $\lambda \in (r-g,r-g+\alpha/\kappa)$. Then,
    \begin{align*}
        \lim_{\lambda \uparrow r-g + \alpha /\kappa} a^* (\lambda) = \overline{a}, \quad \text{and} \quad \lim_{\lambda \uparrow r-g+\alpha /\kappa} b^*(\lambda) = \infty, 
    \end{align*}
    where $\overline{a}$ denotes the threshold characterising the optimal debt issuance policy of the government in equilibrium for large values of $\lambda$, as determined in Theorem \ref{Theorem:NE1}.
\end{corollary}

\section{Comparative Statics Analysis}\label{Sec: Comparative Statics}
In Section 5 we derived the existence and uniqueness of Nash equilibria under different parameter regimes in our model. 
It was revealed that the parameter $\lambda$, measuring the legislative body's time preferences, plays a crucial role on the characterisation of Nash equilibrium. 
While it is optimal not to intervene for large values of $\lambda > r-g+\alpha/\kappa$, the equilibrium strategies are characterised by a pair of thresholds $(a^*,b^*)$ for intermediate values of $\lambda\in (r-g,r-g+\alpha/\kappa)$. 

In this section, we study the sensitivity of these boundaries with respect to some of the model parameters.
In order to highlight the transition from Nash equilibria that are characterised by thresholds $(a^*,b^*)$ (cf.~Theorem \ref{Theorem: Existence and Uniqueness of NE}) to those that prescribe a non-intervention policy for the legislative body (cf.~Theorem \ref{Theorem:NE1}), as stated in Corollary \ref{Coro:ConvNash}, we plot the equilibrium values of $(a^*,b^*)$ as functions of $\lambda$ in the following comparative statics. 

Unless otherwise specified, we fix the following parameter set. \\
\begin{center}
\begin{tabular}{|c|c|c|c|c|c|c|c|c|}
\hline 
    $\rho$ & $\sigma$ & $r$ & $g$ & $\alpha$ & $m$ & $c_1$ & $c_2$ & $\kappa$  \\
    \hline
      $0.3$ & $0.2$ & $0.025$ & $0.02$ & $0.15$ & $0.6$ & $2$ & $1.25$ & $0.6$\\
      \hline
\end{tabular}
\end{center}

\begin{figure}[h]
\centering
\begin{subfigure}{0.35\textwidth}
\includegraphics[height=3.5cm]{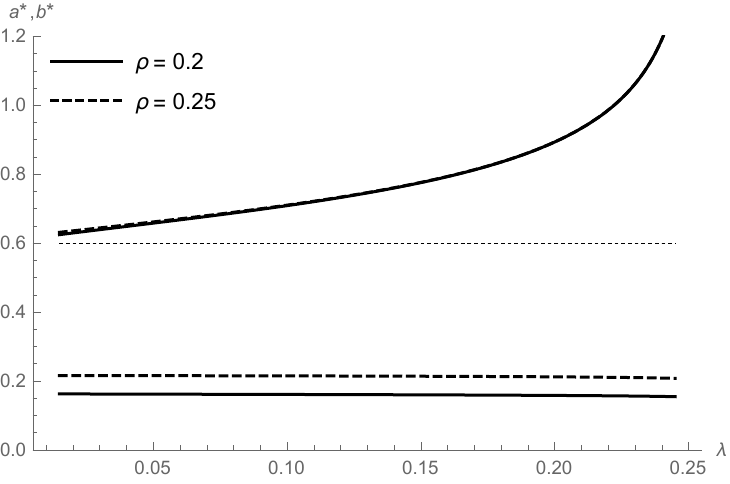} \caption{}
\label{Fig: CompStatrho}
\end{subfigure}\begin{subfigure}{0.35\textwidth}
\includegraphics[height=3.5cm]{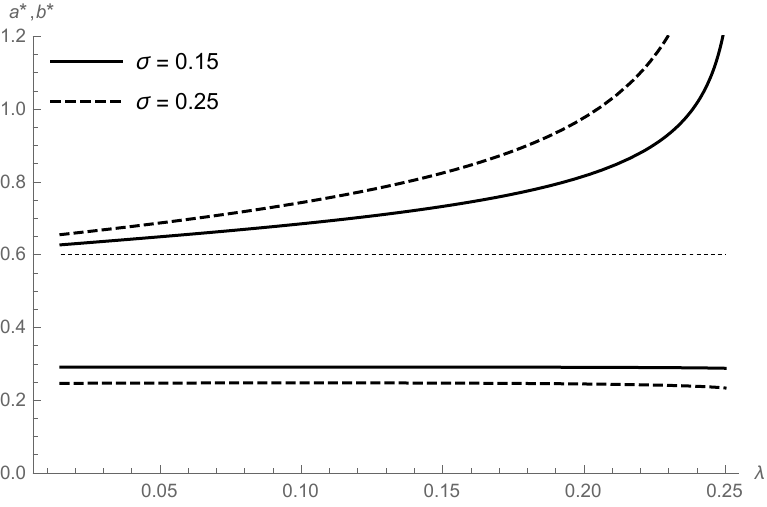}
    \caption{}
\label{Fig: CompStatsigma}
\end{subfigure}

\begin{subfigure}{0.35\textwidth}
\includegraphics[height=3.5cm]{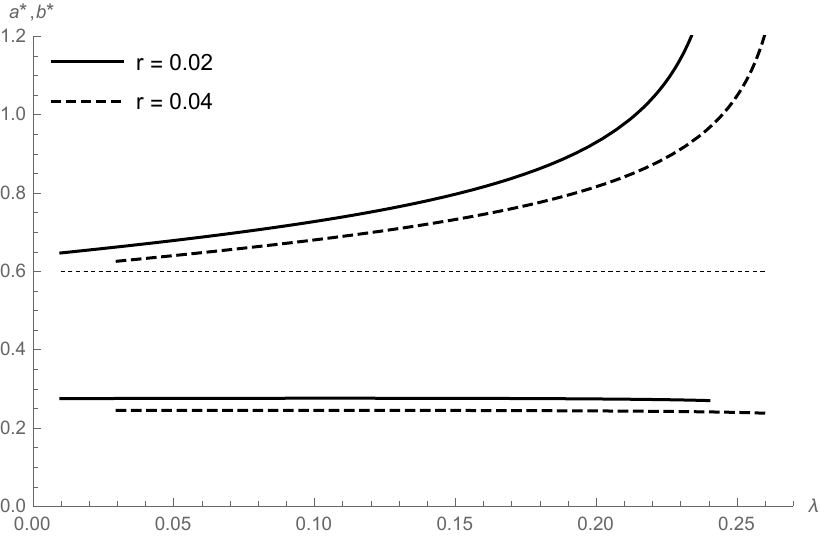}
    \caption{}
\label{Fig: CompStatr}
\end{subfigure}\begin{subfigure}{0.35\textwidth}
\includegraphics[height=3.5cm]{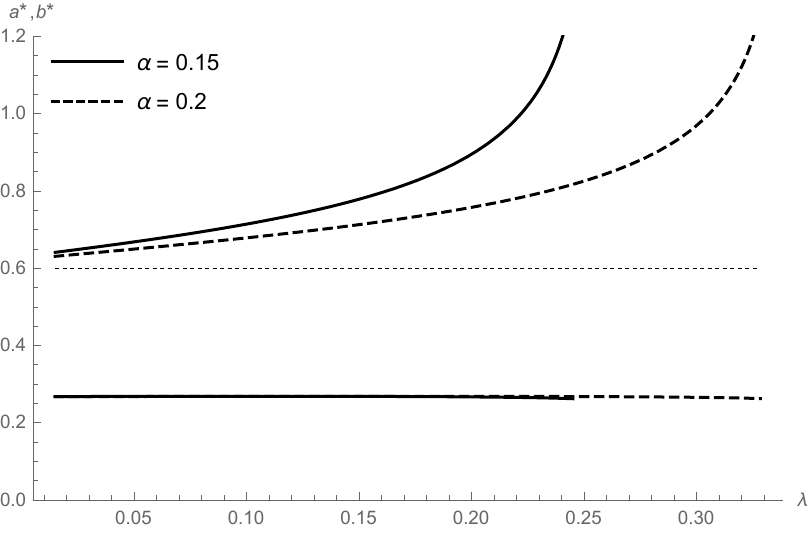}
    \caption{}
\label{Fig: CompStatalpha}
\end{subfigure}
\caption{Sensitivity of the equilibrium values for $a^*$ and $b^*$ (as functions of $\lambda$) with respect to a change in some of the model parameters.}
\label{Fig: CompStatPanel}
\end{figure}

%\subsection
\vspace{3mm}
{\it Sensitivity with respect to $\lambda$.} 
To begin with, it is interesting to study the dependency of the equilibrium values $a^*$ and $b^*$ on the discount factor $\lambda$. 
The numerical sensitivity analysis exhibited in Figure \ref{Fig: CompStatPanel} depicts the optimal intervention thresholds as functions of the legislative body's time preference rate, which here takes values on the interval $\lambda \in (r-g,r-g+\alpha/\kappa)$. The latter guarantees, as shown in the previous analysis, the optimality of a finite debt ceiling mechanism. Some remarks are worth mentioning. Clearly, the equilibrium debt ceiling $b^*$ exhibits a monotonically increasing behaviour as a function of $\lambda$, with the peculiarity of an exploding behaviour $b^* \uparrow + \infty$ for $\lambda \uparrow r-g + \alpha/\kappa$. This illustrates our finding of the smooth transition from a debt ceiling mechanism to a non-intervention policy by the legislative body, as stated in Corollary \ref{Coro:ConvNash}. 
It is interesting to notice that this monotonicity as well as limiting behaviour of $b^*$ does not depend on the other parameters in the model, although some of them influence the interval bounds for the values of $\lambda$, for which the legislative body's optimal strategy is indeed a debt ceiling mechanism. More precisely, increasing (decreasing) the term $r-g$ shifts the interval to the right (left), while the fraction $\alpha/\kappa$ determines the length of the interval $(r-g,r-g+\alpha/\kappa)$. A short discussion on the implications of a shift in these parameters on the optimal strategy of the legislative body (depending on their time preference rate $\lambda$) is given in the subsequent sensitivity study. 
Last, we note that the intervention threshold $a^*$, characterising the optimal debt issuance policy by the government, again illustrates our finding of Corollary \ref{Coro:ConvNash}, in the sense that $a^* \to \overline{a}$ for $\lambda \uparrow r-g+\alpha/\kappa$.

\vspace{3mm}
{\it Sensitivity with respect to governmental time preference rate $\rho$.}
The discount factor $\rho$ serves as a measure on how myopic a government is regarding its debt. Increasing $\rho$ has the effect that the government discounts future costs and revenues more heavily, and thus cares less and less about the future compared to the present. We observe the sensitivity regarding a change in $\rho$ in Figure \ref{Fig: CompStatrho}.
Clearly, the government aims at increasing its debt ratio earlier, at a higher debt-issuance threshold $a^*$, whenever its subjective discount rate increases. 
The legislative body reacts to such an increase by increasing the debt ceiling as well, although we observe Figure \ref{Fig: CompStatrho} that the equilibrium value $b^*$ is relatively robust with respect to  a change in $\rho$.  

%\subsection
\vspace{3mm}
{\it Sensitivity with respect to debt ratio volatility $\sigma$.}
Increasing volatility increases the fluctuations of the debt ratio. The government and the legislative body adapt by acting on the debt ratio later, which is achieved by the government decreasing its optimal debt-issuance threshold and by the legislative body increasing the debt ceiling. %its threshold. 
We can observe this in Figure \ref{Fig: CompStatsigma}.

%\subsection
\vspace{3mm}
{\it Sensitivity with respect to the interest rate $r$ on government debt.}
Increasing interest rates on public debt result in holding debt getting more costly for the government, which in turn increases the drift of the debt ratio. 
Clearly, it is optimal for the government to increase its debt at a later stage, which is achieved by decreasing its debt-issuance %reflection 
threshold, as observed in Figure \ref{Fig: CompStatr}.

In the equilibrium, the legislative body also decreases the debt ceiling, %its threshold, 
since countries with a higher cost of debt are more in danger of defaulting. 
Contrary, if interest rates decrease, the legislative body can be more flexible, since it is optimal to increase the debt ceiling. Intuitively, in such a case, the growth of GDP helps containing the debt ratio without interventions. 
Furthermore, we note that an increase in the interest rate $r$ shifts the interval of $\lambda$-values, for which the optimal strategy of the legislative body prescribes to set a finite debt ceiling, to the right. It follows that a legislative body with a fixed time preference rate could change its optimal strategy from a non-intervention policy to a debt ceiling mechanism, if the government's interest rate on debt increases. Notice that an increase in the country's GDP growth rate has the contrary effect, thus implying that fast growing economies could allow a larger deficit without interventions from a legislative body.

%\subsection
\vspace{3mm}
{\it Sensitivity with respect to the tax compliance factor $\alpha$.} 
The parameter $\alpha$ denotes the tax compliance factor, which measures the willingness to pay taxes within a country. This factor can be chosen freely, hence, the legislative body can account for the fact that some countries have a low probability of default, even though holding a lot of debt. In Figure \ref{Fig: CompStatalpha} we observe the sensitivity of the equilibrium values $a^*$ and $b^*$ with respect to a change in $\alpha$. Clearly, if the tax compliance factor increases (which implies decreasing willingness to pay tax), the legislative body faces stronger social and political pressure to act via the implementation of a debt ceiling mechanism. This has the consequence that a larger factor $\alpha$ (i) causes the legislative body to act earlier on the debt ratio (by decreasing the debt ceiling $b^*$) and (ii) enlarges the interval of time preference values $\lambda \in (r-g,r-g+\alpha/\kappa)$ for which the legislative body's optimal strategy is to impose a debt ceiling mechanism. The latter implies that, for a fixed time preference $\lambda$, a change in the factor $\alpha$ could incentivise the legislative body to switch from a laissez-faire policy to implementing a debt ceiling. On the other hand, if the assigned likelihood of the country's default decreases (in terms of a decrease in the parameter $\alpha$), the legislative body is willing to postpone interventions by either implementing a larger debt ceiling $b^*$ or even choosing a non-intervention policy.

\section{Conclusions}\label{Sec:Conclusions}

In this paper, we study a model of optimal debt management that captures the strategic interaction between a government and a legislative body aiming to optimally control the debt-to-GDP ratio (also called ``debt ratio") of a country. 
While the government is able to increase the debt ratio in favour of public spending, the legislative body aims to restrict it by potentially installing a debt ceiling mechanism. 
The fact that each player has their own objective and optimises their policy while keeping in mind the other player's actions, results in a two-player non-zero-sum game of singular control. We succeed in proving the existence and uniqueness of a Nash equilibrium within the class $\mathcal{M}$, i.e.~in which at least one player acts according to a Skorokhod reflection type policy. This is achieved by studying the variational inequalities associated to the cost minimization problems of each of the two players, and finally determining a ``fixed point" in their best response maps. 

From a mathematical point of view, given the limited amount of solvable non-zero-sum games of singular control, we believe that our detailed study nicely complements the existing literature. 
Furthermore, our modelling framework, that includes a non-differentiable running cost function in the legislative body's performance criterion and different time preferences for each player, destroys the link with a non-zero-sum game of optimal stopping and thus requires an alternative direct approach when searching for Nash Equilibria. 
Our analysis reveals that one can expect two qualitatively different Nash equilibria in the model, depending on the model's parameters. 
We find that while the Nash equilibrium always prescribes the government to install a debt issuance policy (via a Skorokhod reflection at a threshold $a$), the optimality of a debt-ceiling mechanism relies on the legislative body's time preference rate $\lambda$. 
Notably, for large values of $\lambda$, it is optimal to follow a no-intervention policy. 

In economic terms, the latter can be interpreted as the legislative body becoming more myopic regarding future costs of holding debt and thus the government's probability of default. 
This leads them to refrain from taking any actions, leaving the debt ratio unbounded from above. 
On the other hand, for lower values of $\lambda$, we show that it is indeed optimal to install a debt ceiling mechanism by reflecting the debt ratio downwards at a threshold $b$. 
Hence, in this situation, the debt ratio is kept within an interval with minimal effort by the actions of the two opposing players. 
The sensitivity of these equilibrium thresholds to the model parameters is then studied through a comparative statics analysis. 
Amongst other observations, this reveals that a stricter regulation by the legislative body is optimal for higher values of interest rates or tax compliance factors. 
It is also interesting to observe that for particular changes in the underlying parameters, such as the interest rate, GDP growth rate or tax compliance factor, the legislative body may switch from a debt ceiling mechanism to a laissez-faire policy -- and vice versa.  

We finally note that in this paper, we consider Nash equilibria, but it might be interesting to consider other types of equilibria (non-threshold-type, Stackelberg, Pareto, etc.) as well. 
For example, in the environment of Section \ref{Sec: Government Case II}, we could consider the case where the legislative body is credible enough to commit to keep the debt ratio below a debt ceiling $b$ by implementing the reflection policy $\eta^b$. 
Leaving aside the possibility of renegotiation, this corresponds to the debt ceiling policy of the US Congress. 
The best debt issuance policy is thus to reflect the debt ratio process at the level $a(b)$. 
Then, the legislative body chooses a best debt ceiling $b$ by solving the minimisation problem (cf.~\eqref{Objective: ECB New})
$$
\hat V_2(x) = \inf_b \mathcal{I}_{x, \xi^{a(b)}} (\eta^b).
$$
If there is an optimal $\hat b$ for this problem, the pair $(\hat b,a(\hat b))$ is a Stackelberg equilibrium for our non-zero-sum game. 
We leave the general study of Stackelberg equilibria to future research, but make the following observation. 
Assume the legislative body commits to $b^*$ given in Theorem \ref{Theorem: Existence and Uniqueness of NE}, the best debt issuance policy will be $a^*=a(b^*)$ because $(a^*,b^*)$ is a Nash equilibrium. 
Therefore, the cost $\hat V_2(x)$  satisfies
$$
\hat V_2(x) 
\leq \mathcal{I}_{x, \xi^{a(b^*)}} (\eta^{b^*}) 
= \mathcal{I}_{x,\xi^{a^*}} (\eta^{b^*})
= V_2(x;a^*), 
$$
where the latter is the cost for the legislative body along the unique Nash equilibrium (cf.\ Definition \ref{Def: Nash Equilibrium}). 
In other words, it is in the legislative body's interest to commit to a debt ceiling mechanism.

\appendix
\section{Proofs of results in Section \ref{Section: Governments problem}}\label{Appendix:ProofsGovernment}

\subsection{Proof of Theorem \ref{VT:Gcase1}}\label{Appendix:ProofVTG1} 
We derive the result in a number of steps. 

\vspace{1mm}
{\it Step 1.}
We begin with solving the free-boundary problem \eqref{Def: Boundary problem gov case 1}, by constructing a solution to the ordinary differential equation and imposing the boundary conditions, to obtain a \textit{candidate} value function 
\begin{align}\label{Def: Candidate Ubar1}
\overline{U}_1 (x) := 
\begin{cases}
\overline{U}_1 (\overline{a}) - c_2 (\overline{a} -x), & 0 < x \leq \overline{a} , \\
\overline{D}_1 (\overline{a}) x^{\delta_2} + \frac{1}{2(\rho - 2(r-g)-\sigma^2)} x^2 , & \overline{a} < x ,
\end{cases}
\end{align} 
with $\overline{D}_1 (a)$ and $\overline{a}$ as in \eqref{Def: Dbar1 and a*1}. 
%As a preparation for the proof of Theorem \ref{VT:Gcase1}, we first need to 

\vspace{1mm}
{\it Step 2.}
We then aim at verifying that the function $\overline{U}_1$ of \eqref{Def: Candidate Ubar1} solves the free-boundary problem \eqref{Def: Boundary problem gov case 1} and satisfies the HJB equation \eqref{Def: HJB Ubar1}.
Notice that, in view of the construction of $\overline{U}_1$, it remains to check whether 
$$ 
(i) \; (\mathcal{L}-\rho)\overline{U}_1(x) \geq - \frac12 x^2 \quad \text{for $x \in (0,\overline{a})$} 
\qquad \text{and} \qquad 
(ii) \; \overline{U}_1 '(x) \geq c_2 \quad \text{for $x \geq \overline{a}$.}
$$

{\it Proof of (i)}.
We firstly notice that by construction $(\mathcal{L}-\rho)\overline{U}_1(\overline{a}) = - \frac12 \overline{a}^2$. 
We then observe that $x \mapsto (\mathcal{L} - \rho)\overline{U}_1(x)$ decreases with slope $-c_2(\rho - (r-g))$, while $x \mapsto - \frac12 x^2$ decreases with slope $-x$. 
To conclude $(i)$, it is thus sufficient to prove that $c_2 (\rho - (r-g)) > x$ for all $x \in (0,\overline{a})$, or equivalently that $c_2 (\rho - (r-g)) > \overline{a}$. 
The latter follows straightforwardly from the definition \eqref{Def: Dbar1 and a*1} of $\overline{a}$ and $\delta_2 < 0$. 

{\it Proof of (ii)}.
We then show that $x \mapsto \overline{U}_1 '(x)$ is increasing for $x \geq \overline{a}$, by computing  
\begin{align*}
\overline{U}_1''(x) = \frac{1}{\rho - 2(r-g) -\sigma^2} \Big( 1 - \Big(\frac{x}{\overline{a}}\Big)^{\delta_2-1} \Big) > 0,
\end{align*}
where the latter inequality follows from %Assumption \ref{Assumption: c1 > c2}.(ii) with $p=2$ given by 
\eqref{Assumption:Rho:p=2}. 
Hence, given that $\overline{U}_1'(\overline{a}) = 
c_2$ by construction, we conclude that $(ii)$ holds true. 

\vspace{1mm}
{\it Step 3.}
Finally, we must verify that the obtained solution $\overline{U}_1$ of the HJB equation \eqref{Def: HJB Ubar1} identifies with the value function $\overline{V}_1$ of \eqref{Value function Vbar1(x)}. 
The proof is similar to the one of Theorem \ref{VT:Gcase2} given in Appendix \ref{Appendix:ProofVTG2} ({\it Steps 2--4}), thus it is omitted for brevity.

%%%%%%%%%%%%%%%%%%%%%%%%%%%%%%%%%%%%%%%%%%%%%%%%%%%%%%%%%%%%%%%%%%%%%%%%%%%%%%%%%%%%%%%%%%%%%%%%%%%%%%%%%%%%%%%%%%%%%%%%%%%%%%%%%%%%%%%%%%%%%%%%%%%%%%%%%%%%%%%%%%%%%%%%%%%%%%%%%%%%%%%%%%%%%%%%%%%%%%

\subsection{Proof of Lemma \ref{Lemma: a(b) unique}}\label{Appendix:ProofLemmaBounda(b)} 
We prove each part separately.

\vspace{1mm}
\textit{Proof of part} (i). Regarding the existence and uniqueness of a solution $a(b) \in (0,b)$ solving \eqref{Equ: F = 0}, we straightforwardly calculate
\begin{align*}
\lim_{a \downarrow 0} F(a,b) = -\infty \quad \text{and} \quad F(b,b) = (\delta_1 - \delta_2)(c_1 - c_2)(\rho-2(r-g)-\sigma^2)>0,
\end{align*}
where the latter inequality follows from %Assumption \ref{Assumption: c1 > c2}.(ii) with $p=2$ given by 
\eqref{Assumption:Rho:p=2}. 
Then, the first derivative of $F(a , b)$ with respect to $a$ is given by 
\begin{align*}
\frac{\partial}{\partial a} F(a,b) = \, &a^{-1} \Big[ \Big(\frac{b}{a}\Big)^{\delta_1 -1} - \Big( \frac{b}{a} \Big)^{\delta_2 -1} \Big] %\\& \times 
\big[ c_2 (\delta_1-1)(1-\delta_2)(\rho-2(r-g)-\sigma^2)-(\delta_1-2)(2-\delta_2)a\big],
\end{align*}
which implies
\begin{align}\label{Equ: Deriv. of F wrt a}
\frac{\partial}{\partial a} F(a,b) = 
\begin{cases}
> 0, & %\text{for } 
0<a<\tilde{a} \wedge b,\\
< 0, & %\text{for } %a> \tilde{a},
\tilde{a} \wedge b < a < b,
\end{cases}
\quad \text{and} \quad 
\frac{\partial}{\partial a} F(b,b) = 0 ,
\end{align}
for $\tilde{a}$ defined in \eqref{Def: atilde}; note that, the positivity of $\tilde{a}$ follows from \eqref{Assumption:Rho:p=2}.
As a by-product from the above, $F$ crosses zero only once on $(0,b)$ and we can further conclude that $0 < a(b) < %\leq 
\tilde{a} \wedge b$ and $\frac{\partial}{\partial a} F(a(b),b) > 0$.

\vspace{1mm}
\textit{Proof of part} (ii). Regarding the monotonicity results for $a(\cdot)$, we first derive the following partial derivatives
\begin{align}\label{Equ: b partialb F}
b \frac{\partial}{\partial b} F(a,b) &= (\delta_1 -1) [(2- \delta_2)a - c_2 (1 - \delta_2)(\rho - 2 (r-g)-\sigma^2)]\Big(\frac{b}{a}\Big)^{\delta_1 -1}  \\
& + (\delta_2 -1) [(\delta_1 -2)a -c_2(\delta_1 -1)(\rho - 2(r-g)-\sigma^2)]\Big(\frac{b}{a} \Big)^{\delta_2 -1} - (\delta_1 - \delta_2)b \nonumber
\end{align}
and
\begin{align}\label{Equ: bb partialbb F}
b^2 \frac{\partial^2}{\partial b^2} F(a,b) &= (\delta_1 -1)(\delta_1-2) \big[ (2-\delta_2)a - c_2(\rho - 2(r-g)-\sigma^2)(1- \delta_2) \big] \Big( \frac{b}{a} \Big)^{\delta_1-1} \\
&+ (1-\delta_2)(2-\delta_2) \big[ (\delta_1-2)a - c_2(\rho - 2(r-g) - \sigma^2) (\delta_1-1) \big] \Big( \frac{b}{a}\Big)^{\delta_2-1}. \nonumber
\end{align}
Furthermore, given that $\frac{\partial}{\partial a} F(a(b),b) > 0$ at the value $a(b)$ which satisfies \eqref{Equ: F = 0} due to %Lemma \ref{Lemma: a(b) unique}
part (i), we can obtain the monotonicity of $a(b)$ on $(0,\infty)$ through
\begin{align}\label{Equ: a'(b)}
a'(b) = - \frac{\frac{\partial}{\partial b} F(a(b),b)}{\frac{\partial}{\partial a} F(a(b),b)} \geq 0
\quad \Leftrightarrow \quad 
\frac{\partial}{\partial b} F(a(b),b) \leq 0. 
\end{align}

% \begin{proposition}\label{Lemma: Mono of a(b)}
% Let $a(b)$ be the unique solution to \eqref{Equ: F = 0} as in Lemma \ref{Lemma: a(b) unique}. 
% Then, %$b \mapsto a(b)$ %is first increasing and then decreasing. 
% \begin{equation}\label{a'0-inf}
% \text{$a(\cdot)$ is} 
% \begin{cases} 
% \text{increasing on $(0, \hat b)$}, \\
% \text{decreasing on $(\hat b, \infty)$}, 
% \end{cases}
% \quad \text{where $\hat{b}$ is the unique solution to $\frac{\partial}{\partial b} F (a(\hat{b}),\hat{b})=0$.}
% \end{equation}
% Moreover, $b \mapsto a(b)$ is concave on the interval $(0,\hat{b})$ and
% $\lim_{b\to 0} a(b) = 0$ as well as $\lim_{b\to \infty} a(b) = \overline{a}$, where $\overline{a}$ is the optimal debt-issuance threshold defined by \eqref{Def: Dbar1 and a*} in Case (I) of non-intervention by a legislative body. 
% \end{proposition}
% \begin{proof} 
In the following, we fix $b\in \mathbb{R}_+$ and consider the uniquely defined $a(b) \in (0,\tilde{a} \wedge b)$ given by the solution to \eqref{Equ: F = 0}. We distinguish two cases depending on the location of the fixed $a(b)$ relative to $\overline{a}$ defined in \eqref{Def: Dbar1 and a*1}. 

{\it Case} (a): $a(b) \leq \overline{a}$. 
Given that $\overline{a} < \tilde{a}$ due to the assumption in \eqref{Assumption:Rho:p=2} and the definitions \eqref{Def: Dbar1 and a*1} of $\overline{a}$ and \eqref{Def: atilde} of $\tilde{a}$, we observe that 
$$
\frac{\partial^2}{\partial b^2} F (a(b),x) < 0, \quad \text{for all $x > a(b)$} 
\quad \Rightarrow \quad 
x \mapsto \frac{\partial}{\partial b} F (a(b),x) \text{ is strictly decreasing on } (a(b),\infty). $$ 
Since $\frac{\partial}{\partial b} F(a(b),a(b)) = 0$ due to \eqref{Equ: b partialb F}, we have $\frac{\partial}{\partial b} F (a(b),x) < 0$ for all $x > a(b)$. 
Combining this with the fact that $b > a(b)$, %for all $b \in \mathbb{R}_+$ due to Lemma \ref{Lemma: a(b) unique}, 
we conclude from \eqref{Equ: a'(b)} that $a'(b) > 0$ for all $b \in \mathbb{R}_+$ s.t. $a(b) \leq \overline{a}$.
This yields that 
\begin{equation}\label{a'0-barb}
\text{$a(\cdot)$ is increasing on $(0,\overline{b}]$, \qquad where $\overline{b}$ is such that $a(\overline{b})=\overline{a}$.}
\end{equation}
Also, we observe that $a'(\overline{b})>0$ and $a(b)>\overline{a}$ for all $b>\overline{b}$.

{\it Case} (b): $a(b) > \overline{a}$. 
We firstly note from Case (a) that this is realised when $b>\overline{b}$.
We observe that $\frac{\partial^2}{\partial b^2} F (a(b),x) \leq 0$ if and only if $x \leq \tilde{x}$, where 
\begin{align*}
\tilde{x} := \Big( \frac{(2 - \delta_2)(1 - \delta_2)(c_2 (\delta_1-1)(\rho - 2(r-g)-\sigma^2) - (\delta_1-2)a(b)) a(b)^{\delta_1 - \delta_2}}{(\delta_1-2)(\delta_1 -1) ((2-\delta_2)a(b) - c_2 (1 - \delta_2)(\rho - 2(r-g)-\sigma^2)} \Big)^{\frac{1}{\delta1 - \delta2}} > a(b), 
\end{align*}
which is well-defined since $a(b) \in (\overline{a} , \tilde{a})$. %due to Lemma \ref{Lemma: a(b) unique}. 
To show the inequality via contradiction, assume that $\tilde{x} \leq a(b)$. 
Then, $\frac{\partial^2}{\partial b^2} F (a(b),x) \geq 0$ and hence $x \mapsto \frac{\partial}{\partial b} F(a(b),x)$ is increasing for all $x \geq a(b)$. 
But since $\frac{\partial}{\partial b} F(a(b),a(b)) = 0$, it would follow that $x \mapsto F(a(b),x)$ is increasing on $(a(b),\infty)$, which is a contradiction to $F(a(b),b) = 0$, given that $a(b)<b$.  
Therefore, 
\begin{align*}
\frac{\partial^2}{\partial b^2} F (a(b),x) 
\begin{cases}
\leq 0 , & %\text{for } 
a(b) \leq x \leq \tilde{x}, \\
\geq 0, & %\text{for } 
x \geq \tilde{x}
\end{cases}
\quad \Rightarrow \quad 
x \mapsto \frac{\partial}{\partial b} F (a(b),x) \text{ is}
\begin{cases} 
\text{decreasing on } (a(b),\tilde{x}), \\
\text{increasing on } (\tilde{x},\infty). 
\end{cases}
\end{align*}
Combining this with the fact that $\frac{\partial}{\partial b} F (a(b),a(b)) = 0$ and $\lim_{x \to \infty} \frac{\partial}{\partial b} F (a(b),x) = + \infty$, we conclude that $\frac{\partial}{\partial b} F (a(b),x) = 0$ admits a unique solution on $(a(b), \infty)$, denoted by $x_m(b) \in (\tilde{x}, \infty)$. 
Hence, 
\begin{align*}
x \mapsto \frac{\partial}{\partial b} F (a(b),x)  
\begin{cases} 
< 0 , %\text{ for } 
&x \in (a(b),x_m(b)), \\
> 0 , %\text{ for } 
&x \in (x_m(b),\infty) 
\end{cases} 
\quad \Rightarrow \quad 
x \mapsto F (a(b),x) \text{ is}
\begin{cases} 
\text{decreasing on } (a(b),x_m(b)), \\
\text{increasing on } (x_m(b),\infty). 
\end{cases}
\end{align*}
Given that $F(a(b),a(b)) > 0$ (see the first part of the proof) and $\lim_{x \to \infty} F(a(b),x) = +\infty$, we conclude that there exist at most two solutions to $F(a(b),x) = 0$ 
and due to \eqref{Equ: a'(b)} that $a'(b)$ changes sign %only 
once. 
This implies -- in view of the conclusion $a'(\overline{b})>0$ in Case (a) -- that $a(\cdot)$ is either increasing on the whole $(\overline{b}, \infty)$, or it is increasing only on $(\overline{b}, \hat b)$ and then decreasing on $(\hat b, \infty)$, where $\hat b \in (\overline b, \infty)$ would be satisfying 
$$a'(\hat b)=0 
\quad \Leftrightarrow \quad 
\frac{\partial}{\partial b} F (a(\hat b),\hat b)=0
\quad \Leftrightarrow \quad 
\hat b = x_m(\hat b).
$$
In order to show that such a $\hat b$ always exists, we study the system of equations $F(\hat{a},\hat{b}) = 0 = \frac{\partial}{\partial b} F (\hat{a},\hat{b})$ (cf.\ equations above), which is equivalent to
\begin{align}\label{Equ: Coupled for ahat bhat}
\begin{cases}
J_{1,2}(\hat{a}) = J_{1,1}(\hat{b}), \\
J_{2,2} (\hat{a}) = J_{2,1}(\hat{b}),
\end{cases}
\hspace*{1cm} \text{where} ~~ J_{i,j} (x) := \frac{(\delta_i -2)x - c_j (\delta_i-1)(\rho-2(r-g)-\sigma^2)}{x^{\delta_{3-i}-1}}
\end{align}
It can be shown (see, e.g. \cite{FeRo}) that the system \eqref{Equ: Coupled for ahat bhat} admits a unique solution $(\hat{a},\hat{b}) \in \mathbb{R}_+^2$, where 
\begin{equation}\label{a'barb-inf}
\text{$\hat{a}=a(\hat b)$, such that $a'(\hat b)=0$ and }
\text{$a(\cdot)$ is} 
\begin{cases} 
\text{increasing on $(\overline{b}, \hat b)$}, \\
\text{decreasing on $(\hat b, \infty)$}. 
\end{cases} 
\end{equation}

The monotonicity then follows by combining \eqref{a'0-barb} and \eqref{a'barb-inf}. 
Furthermore, this monotonicity together with the fact that for every choice of $b \in \mathbb{R}_+$, there always exists a best response $a(b)$ that satisfies \eqref{Equ: F = 0}--\eqref{Def: F(a,b)}, due to %Lemma \ref{Lemma: a(b) unique}
part (i), then yields $\lim_{b\to\infty} a(b) = \overline{a}$.

\vspace{1mm}
\textit{Proof of part} (iii).
Regarding the concavity of $a(\cdot)$ on the interval $(0,\hat{b})$ we investigate the term 
\begin{align*} %\label{Def: a''}
    a''(b) = \frac{2 F_a (a(b),b) F_b (a(b),b) F_{ab}(a(b),b) -  F_b^2(a(b),b) F_{aa} (a(b),b) - F_a^2 (a(b),b) F_{bb}(a(b),b) }{F_a^3 (a(b),b)}, 
\end{align*}
where we set $F_a (a,b) := \partial_a F(a,b)$ (and the other terms analogously). 
We first notice that, due to \eqref{Equ: Deriv. of F wrt a}, we have $F_a(a(b),b)>0$.
Direct computation yields
\begin{align*} 
&2 F_a (a(b),b) F_b (a(b),b) F_{ab}(a(b),b) -  F_b^2(a(b),b) F_{aa} (a(b),b) - F_a^2 (a(b),b) F_{bb}(a(b),b) \nonumber \\
&= \Big\{ \frac{1}{a(b) b} (2-\delta_2) \Big( \frac{b}{a(b)}\Big)^{\delta_1-1} + \frac{1}{a(b) b} (\delta_1-2)\Big(\frac{b}{a(b)}\Big)^{\delta_2-1} - (\delta_1-\delta_2)\frac{1}{a(b)^2} \Big\} \nonumber\\
&\quad\times \Big\{ (\delta_1-\delta_2) (\rho-2(r-g)-\sigma^2) \big[ c_2 (\delta_1-1)(1-\delta_2) (\rho-2(r-g)-\sigma^2) - (\delta_1-2)(2-\delta_2)a(b) \big]  \nonumber \\
&\quad\qquad \times \frac{1}{b} \Big[ c_2(\delta_1-\delta_2) \Big( \frac{b}{a(b)}\Big)^{\delta_1 + \delta_2-2} %\Big(\frac{b}{a}\Big)^{\delta_2-1}
- c_1 (\delta_1-1) \Big(\frac{b}{a(b)}\Big)^{\delta_1-1} - c_1(1-\delta_2)\Big(\frac{b}{a(b)}\Big)^{\delta_2-1} \Big] \nonumber \\
&\quad\qquad + c_2(1-\delta_2)(\delta_1-1)(\rho-2(r-g)-\sigma^2) \Big[ \Big(\frac{b}{a(b)}\Big)^{\delta_1-1} - \Big( \frac{b}{a(b)}\Big)^{\delta_2-1} \Big] F_b (a(b),b) \Big\}.\nonumber
\end{align*}
Some straightforward calculations reveal that the first term on the above right-hand side (second line) is strictly positive, 
while the term in the third line is strictly positive due to $a(b) \leq \tilde{a}$. 
Moreover, the term in the fourth line is strictly negative, which easily follows upon using $a(b) \leq b$. Finally, we notice that for $b \leq \hat{b}$, we have $F_b (a(b),b) < 0$ (see \eqref{a'barb-inf}). Combining these facts, we conclude that indeed $a''(b) < 0$ for $b\in (0,\hat{b})$ and the claim follows.  \qed
%\end{proof}

%%%%%%%%%%%%%%%%%%%%%%%%%%%%%%%%%%%%%%%%%%%%%%%%%%%%%%%%%%%%%%%%%%%%%%%%%%%%%%%%%%%%%%%%%%%%%%%%%%%%%%%%%%%%%%%%%%%%%%%%%%%%%%%%%%%%%%%%%%%%%%%%%%%%%%%%%%%%%%%%%%%%%%%%%%%%%%%%%%%%%%%%%%%%%%%%%%%%%%

\subsection{Proof of Theorem \ref{VT:Gcase2}}\label{Appendix:ProofVTG2}
We derive the result in a number of steps. 
In particular, we show in {\it Step 1} that the candidate value function \eqref{Def: Cand VF V1}, with $a(b)$ solving \eqref{Equ: F = 0} as in Lemma \ref{Lemma: a(b) unique}, indeed solves the HJB equation \eqref{Def: HJB U1(x,b)}, and in {\it Steps 2--4} that the latter identifies with value function $V_1$ in \eqref{Value function V1(x,b)}.  

\vspace{1mm}
{\it Step 1.}
By construction, we have $(\mathcal{L} - \rho) U_1 (x;b) = - \frac12 x^2$ for $x \in (a(b),b)$,  $U_1'(x;b) = c_2$ for $x \in (0,a(b))$ and $U_1'(x;b) = c_1 > c_2$ for $x \geq a(b)$. 
Therefore, it remains to show that: 
$$
(i) \; (\mathcal{L} - \rho) U_1 (x;b) \geq - \frac12 x^2 \quad \text{for $x \in (0,a(b))$} 
\qquad \text{and} \qquad 
(ii) \; U_1'(x;b) \geq c_2 \quad \text{for $x \in (a(b),b)$.}
$$
In the following, we fix $b \in \mathbb{R}_+$, so that $a(b)$ is the (fixed) unique solution to \eqref{Equ: F = 0} according to Lemma \ref{Lemma: a(b) unique}. 

{\it Proof of (i)}.
For $x \in (0,a(b))$, we get 
\begin{align*}
(\mathcal{L}- \rho) U_1 (x;b) &= (r-g)x c_2 - \rho U_1(a(b);b) + \rho c_2 (a(b) -x).
\end{align*}
Clearly, $x \mapsto (\mathcal{L} -\rho) U_1 (x;b)$ decreases with slope $- c_2 (\rho - (r-g))$, while $x \mapsto - \frac12 x^2$ decreases with slope $-x$. Since $(\mathcal{L} - \rho)U_1 (a(b);b) = - \frac12 a(b)^2$, it is sufficient to then show that $c_2 (\rho - (r-g)) > x$, for all $x \in (0,a(b))$. The latter is true due to \eqref{Def: atilde}, thus $(i)$ holds true.

{\it Proof of (ii)}.
For $x \in (a(b),b)$, we can calculate 
\begin{align*}
U_1'(x;b) &= \frac{(\delta_2 -2) a(b) - c_2 (\delta_2 -1)(\rho - 2(r-g)-\sigma^2)}{(\delta_1 - \delta_2)(\rho - 2(r-g) - \sigma^2)}\Big( \frac{x}{a(b)}\Big)^{\delta_1 -1} \\
& - \frac{(\delta_1 -2) a(b) - c_2 (\delta_1 -1) (\rho -2(r-g)-\sigma^2)}{(\delta_1 - \delta_2)(\rho - 2(r-g)-\sigma^2)} \Big(\frac{x}{a(b)}\Big)^{\delta_2 -1} + \frac{x}{\rho - 2(r-g) - \sigma^2}.
\end{align*}
Combining this with the definition \eqref{Def: F(a,b)} of $F$, 
we notice that 
$$ U_1'(x;b) \geq c_2 %if and only if
\quad \Leftrightarrow \quad F(a(b),x) \leq (c_1 -c_2) (\delta_1 -\delta_2) (\rho - 2(r-g) - \sigma^2). $$ 
To prove $(ii)$, given that  
\begin{align*}
F(a(b),a(b)) = (c_1 -c_2) (\delta_1 -\delta_2) (\rho - 2(r-g) - \sigma^2) \qquad \text{and} \qquad F(a(b),b) = 0,
\end{align*}
it is sufficient to have that $x \mapsto F(a(b),x)$ first decreases and changes sign at most once %for $x \geq a(b)$. 
in $(a(b), b)$. 
The latter was shown in the proof of Lemma \ref{Lemma: a(b) unique}, thus $(ii)$ holds true. 

\vspace{1mm}
\textit{Step 2.} 
Let $x \in \mathbb{R}_+$ and $\xi \in \mathcal{A}_{\eta^b}$ such that $(\xi , \eta^b ) \in \mathcal{A}$. 
For $n \geq 1$, we let $\tau_n := \inf \{ t \geq 0:~X_t^{0,\eta^b} \geq n\}$. %~\P$-a.s.. 
Since $U_1 \in C^2 (0,b)$, we can apply It\^{o}-Meyer formula to the process $e^{- \rho \tau_n} U_1(X_{\tau_n}^{\xi , \eta^b};b )$ on $[0,\tau_n]$ and obtain
\begin{align}\label{Equ: Ito in VerTh of V1}
e^{-\rho \tau_n} U_1(X_{\tau_n}^{\xi,\eta^b} ; b) - U_1(x;b) 
= &\int_0^{\tau_n} e^{-\rho s} (\mathcal{L} -\rho)U_1(X_s^{\xi , \eta^b} ; b) ds + \sigma \int_0^{\tau_n} e^{-\rho s} U_1'(X_s^{\xi, \eta^b} ; b) d W_s \nonumber \\
&- \int_0^{\tau_n} e^{-\rho s} X_s^{\xi ,\eta^b} U_1 '(X_s^{\xi, \eta^b};b) d \eta_s^{c,b} 
+\int_0^{\tau_n} e^{-\rho s} X_s^{\xi,\eta^b} U_1 '(X_s^{\xi,\eta^b};b )d\xi_s^c \nonumber \\
&+ \sum_{s \leq \tau_n} e^{-\rho s} \big( U_1(X_{s}^{\xi , \eta^b};b) - U_1(X_{s-}^{\xi,\eta^b};b)\big).  
\end{align} 
Clearly, for any $s \in (0,\tau_n]$ we have $0 < X_s^{\xi,\eta^b} < n \wedge b$ and thus continuity of $U_1'$ implies that the second term in \eqref{Equ: Ito in VerTh of V1} is a martingale. Furthermore, since $(\xi, \eta^b) \in \mathcal{A}$, the last term in \eqref{Equ: Ito in VerTh of V1} rewrites as 
\begin{align*}
\sum_{s < \tau_n} e^{- \rho s } \big( U_1(X_{s}^{\xi , \eta^b} ; b)- U_1(X_{s-}^{\xi , \eta^b} ;b) \big) =  \sum_{s < \tau_n} e^{- \rho s } \big( U_1(X_{s}^{\xi , \eta^b} ; b)- U_1(X_{s-}^{\xi , \eta^b} ;b) \big) \big[ \one_{ \{ \Delta \xi_s > 0 \} } + \one_{ \{ \Delta \eta^b_s > 0 \} } \big]
\end{align*}
and 
\begin{align} \label{Equ: Sum for circ u notation} 
\begin{split}
\big( 
U_1(X_{s}^{\xi , \eta^b} ; b) - U_1(X_{s-}^{\xi, \eta^b}& ;b) \big) \one_{ \{ \Delta \xi_s > 0 \} } 
= U_1 (e^{\Delta \xi_s} X_{s-}^{\xi,\eta^b};b) - U_1 (X_{s-}^{\xi,\eta^b};b) \\
&= \int_0^{\Delta \xi_s} \frac{\partial U_1 ( e^u X_{s-}^{\xi,\eta^b}  ; b)}{\partial u} du 
=  \int_0^{ \Delta \xi_s } e^u X_{s-}^{\xi, \eta^b} U_1 ' (e^u X_{s-}^{\xi , \eta^b} ; b) du \\
 \big( U_1(X_{s}^{\xi , \eta^b} ; b) - U_1(X_{s-}^{\xi , \eta^b}& ;b) \big) \one_{ \{ \Delta \eta^b_s > 0 \} } 
 = U_1(e^{-\Delta \eta_s^b} X_{s-}^{\xi , \eta^b}  ; b)- U_1(X_{s-}^{\xi , \eta^b} ;b) \\
 &= \int_0^{\Delta \eta_s^b} \frac{\partial U_1 ( e^{-u} X_{s-}^{\xi,\eta^b} ;b) }{\partial u} d u 
 = -  \int_0^{ \Delta \eta^b_s } e^{-u} X_{s-}^{\xi, \eta^b} U_1 ' (e^{-u} X_{s-}^{\xi , \eta^b} ; b) du
\end{split} 
%\label{Equ: Sum for circ d notation}
\end{align}
Taking expectations, rearranging terms and using the notation introduced in \eqref{Def: Operator circ up} %-\eqref{Def: Operator circ down} 
we can thus write \eqref{Equ: Ito in VerTh of V1} as  
\begin{align}\label{Equ: U1 leq Obj for VT}
 U_1(x;b) &= \E_x \Big[ e^{-\rho \tau_n} U_1(X_{\tau_n}^{\xi,\eta^b};b) -  \int_0^{\tau_n} e^{-\rho s} (\mathcal{L} - \rho) U_1 (X_s^{\xi,\eta^b} ;b) ds + \int_0^{\tau_n} e^{-\rho s} X_s^{\xi,\eta^b} U_1'(X_s^{\xi,\eta^b} ;b) \circ_d d \eta_s^b \nonumber \\ 
& \hspace*{3.95cm} - \int_0^{\tau_n} e^{-\rho s} X_s^{\xi,\eta^b} U_1'(X_s^{\xi,\eta^b};b) \circ_u d\xi_s \Big] \nonumber \\
&\leq \E_x \Big[  e^{-\rho \tau_n} U_1(X_{\tau_n}^{\xi,\eta^b};b)  + \int_0^{\tau_n} e^{-\rho s} h(X_s^{\xi,\eta^b}) ds + c_1 \int_0^{\tau_n} e^{-\rho s} X_s^{\xi,\eta^b} \circ_u d \eta_s^b  \nonumber \\
& \hspace*{3.95cm} -  c_2 \int_0^{\tau_n} e^{-\rho s} X_s^{\xi,\eta^b} \circ_u d\xi_s \Big],
\end{align}
where, in the latter inequality, we exploit the fact that $U_1$ solves the free-boundary problem \eqref{Free Boundary Problem Government}, as stated above. %\nr{[Here and before, we need to be a bit more careful about what $U_1$ solves: HJB vs FB problem (these are not necessarily equivalent) - I will fix this]}
By admissibility of $\xi$ we have that the right-hand side of \eqref{Equ: U1 leq Obj for VT} is finite $\P$-a.s., and we notice that Assumption \ref{Assumption: c1 > c2}.(ii) %{Assumption: Rho} 
guarantees 
\begin{align*}
\lim_{n \to \infty} \E_x \big[ e^{- \rho \tau_n} U_1(X_{\tau_n}^{\xi,\eta^b};b)\big] = 0.
\end{align*}
Then, noticing that $\tau_n \uparrow \infty$, $\P$-a.s., and taking limits in \eqref{Equ: U1 leq Obj for VT}, we can use Assumptions \ref{Assumption: Cost function h}.(iii) and
\ref{Assumption: c1 > c2}.(ii), %\ref{Assumption: Rho} 
employ the dominated convergence theorem, and conclude that 
\begin{align*}
U_1(x;b) \leq \E_x \Big[  \int_0^{\infty} e^{-\rho s} h(X_s^{\xi,\eta^b}) ds + c_1 \int_0^{\infty} e^{-\rho s} X_s^{\xi,\eta^b} \circ_d d \eta_s^b  
- c_2 \int_0^{\infty} e^{-\rho s} X_s^{\xi,\eta^b} \circ_u d\xi_s  \Big].
\end{align*}
Since $\xi \in \mathcal{A}_{\eta^b}$ was arbitrary, we have $U_1(x;b) \leq V_1(x;b)$ on $\mathbb{R}_+$. 

\vspace{1mm}
\textit{Step 3.} 
We can repeat the above arguments from {\it Step 2}, but now fix the control strategy $\xi^{a(b)}$ of \eqref{Def: Optimal Control a(b)}. 
Since $X_t^{\xi^{a(b)} , \eta^b} \in [a(b),b]$, $\P$-a.s., for all $t > 0$, and $U_1(X_t^{\xi^{a(b)},\eta^b} ;b) = c_2$ on  supp$\{d \xi_t^{a(b)} \}$, the inequality in \eqref{Equ: U1 leq Obj for VT} becomes an equality and hence, employing dominated convergence arguments as before, we observe 
\begin{align*}
U_1(x;b) &= \E_x \Big[  \int_0^{\infty} e^{-\rho s} h(X_s^{\xi^{a(b)},\eta^b}) ds + c_1 \int_0^{\infty} e^{-\rho s} X_s^{\xi^{a(b)},\eta^b} \circ_d d \eta_s^b 
- c_2 \int_0^{\infty} e^{-\rho s} X_s^{\xi^{a(b)},\eta^b} \circ_u d\xi_s^{a(b)}  \Big] \\
 &\geq V_1(x;b).
\end{align*}

\textit{Step 4.} 
Combining the results from {\it Steps 2} and {\it 3} then concludes that $U_1(x;b) = V_1(x;b)$ and $\xi^{a(b)}$ of \eqref{Def: Optimal Control a(b)} is an optimal control strategy for problem \eqref{Value function V1(x,b)}.

\section{Proofs of results in Section \ref{Section: The institutions problem}}

\subsection{Proof of Theorem \ref{VT:Icase1}}\label{Appendix:ProofVTI1}
We derive the result in a number of steps. 

\vspace{1mm}
{\it Step 1.}
We begin by solving the fixed-boundary problem \eqref{Def: Free Boundary Problem ECB b infty}. To this end, we construct a solution to the ordinary differential equation an impose the stated boundary conditions. We then obtain a \textit{candidate} value function
\begin{align}\label{Def: Candidate U2 b infty}
\overline{U}_2 (x;a)  &= 
\begin{cases} 
\overline{U}_2 (a;a), & 0 < x \leq a; \\
\overline{D}_2 (a) x^{\theta_2} + H(x),\hspace*{1.7cm} & a<x, 
\end{cases}
\end{align}
where $\overline{D}_2(\cdot)$ and $H(\cdot)$ are as in \eqref{Def:FunctionsforVTI1}. 

\vspace{1mm}
{\it Step 2.}
We then aim at verifying that $\overline{U}_2$ of \eqref{Def: Candidate U2 b infty} solves the free-boundary problem \eqref{Def: Free Boundary Problem ECB b infty} and satisfies the HJB equation \eqref{Def: HJB ECB}.
In view of the construction of $\overline{U}_2$, %the first and third line in \eqref{Def: Free Boundary Problem ECB b infty} are fulfilled. It is thus left 
it remains to check that $0 \leq \frac{\partial}{\partial x} \overline{U}_2(x;a) < \kappa$, %as well as $\frac{\partial}{\partial x} \overline{U}_2^+ < \kappa$ 
for all $x > a$. 

Straightforward calculations lead to
\begin{align*}
%\frac{\partial}{\partial x} \overline{U}_2 (x;a) &= \frac{\alpha}{\lambda - (r-g)} \Big( 1 - \Big( \frac{x}{a}\Big)^{\theta_2 -1} \Big) < \kappa, \\\text{and}\hspace*{1cm} 
\frac{\partial}{\partial x} \overline{U}_2 (x;a) &= \alpha \int_0^\infty e^{-(\lambda - (r-g))t} \Big( \Phi (d_1 (x,t)) - \Big( \frac{x}{a}\Big)^{\theta_2 -1} \Phi (d_1 (a,t)) \Big) dt < \kappa,
\end{align*}
where the inequality follows from %Assumption \ref{Assumption: c1 > c2}.(iii), %\ref{Assumption: lambda}, and 
the facts that $x > a$ and $\lambda > r-g + \frac{\alpha}{\kappa}> r-g$ under Case (\textrm{I}). 
Furthermore, Assumption \ref{Assumption: c1 > c2}.(iii) %Assumption \ref{Assumption: lambda} 
guarantees that
$\frac{\partial}{\partial x} \overline{U}_2(x;a) \geq 0$, 
%as well as $\frac{\partial}{\partial x} \overline{U}_2^{+}(x;a) \geq 0$
for all $x >a$.

\vspace{1mm}
{\it Step 3.}
Finally, we must verify that indeed ${V}_2 = \overline{U}_2$ and that not intervening is an optimal debt management strategy. 
The proof follows the lines of the proof of Theorem \ref{VT:Icase2} ({\it Steps 2--4}), and is thus omitted for brevity.

\subsection{Proof of Lemma \ref{Lem:b(a)}}\label{Appendix:ProofLemmaBoundb(a)}
We prove each part separately. 

%\subsubsection*{Existence and Uniqueness}
\vspace{1mm}
{\it Proof of part }(i). 
Regarding the existence and uniqueness of a solution $b(a) \in (a,\infty)$ solving \eqref{Equ: G=0}, we first conclude from the representation of $G$ in \eqref{Def: G(a,b)} that 
\begin{align*}
%G(a,a) &= (\theta_1- \theta_2)\kappa (\lambda - (r-g)) > 0,  \qquad \lim_{b \to \infty} G(a,b) = - \infty, \\
%\text{and} \hspace*{1cm} 
G(a,a) &= k > 0, \qquad \text{and} \qquad \lim_{b \to \infty} G(a,b) = - \infty,
\end{align*}
where the latter follows precisely from the fact that $\lambda < r-g + \frac{\alpha}{\kappa}$ in Case (\textrm{II}). Hence, there exists a solution to the equation \eqref{Equ: G=0}. 
%Moreover, straightforward calculations reveal 
Moreover, it follows from the following expression 
\begin{align}\label{Equ: G partial b}
\frac{\partial}{\partial b} G (a,b) &= \Big[ \Big( \frac{b}{a} \Big)^{1-\theta_2} - \Big(\frac{b}{a}\Big)^{1-\theta_1} \Big]  \frac{(\theta_1-1)(1-\theta_2)}{(\theta_1 - \theta_2)(\lambda - (r-g)) b} \big[ \kappa (\lambda - (r-g)) - \alpha \one_{ \{b>m\} } \big] ,
\end{align}
that 

\begin{align}\label{Equ:GPartb}
\frac{\partial}{\partial b} G (a,b) 
&=
\begin{cases}
> 0, & a<b<a \vee m, \\ %b<m, \\
< 0, & a \vee m < b. %b>m,
\end{cases} 
\end{align}
Therefore, for any $a \in \mathbb{R}_+$, there exists a unique $b(a) \in (a \vee m, \infty)$ %\in (a,\infty)$ 
such that $G(a,b(a)) = 0$. 
Moreover, due to \eqref{Equ: G partial b}, we conclude that $\frac{\partial}{\partial b} G(a, b(a))<0$.

\vspace{1mm}
\textit{Proof of part} (ii). Let $b_0$ be defined as in \eqref{Def: b0} and fix $a \in \mathbb{R}_+$. In order to show \eqref{Def: b0}, we examine two cases of $a$-values. 

{\it Case} (a): $a \geq b_0$. In this case, we immediately have $b(a)>b_0$. 

{\it Case} (b): $a<b_0$. 
Notice that simple comparison arguments yield that $b_0 > m$ and it can be shown that $G(a,b_0) \geq 0$. 
We then assume (aiming for contradiction) that $b(a) \in (a \vee m, b_0)$. 
%Due to \eqref{Equ: G partial b}, it is clear that $b(a) >m$ for any $a\in \mathbb{R}_+$, and simple comparison yields $b_0 > m$. Let $a \in \mathbb{R}_+$ and assume that $b(a) \in (m,b_0)$. 
Since $b \mapsto G(a,b)$ is strictly decreasing for $b>a \vee m$ (see \eqref{Equ:GPartb}), it follows that
\begin{align*}
G(a,b(a)) &> G(a,b_0) \geq 0,
\end{align*}
which is a contradiction, and it thus follows that $b(a) \geq b_0$. 

\vspace{1mm}
\textit{Proof of part} (iii). We are now in position to obtain the monotonicity of $b(\cdot)$ on $(0,\infty)$. 
Given that $\frac{\partial}{\partial b} G(a,b(a)) < 0$ at the value $b(a)$ which satisfies \eqref{Equ: G=0}, due to the reasoning above, we have that
\begin{align}\label{Def: Derivative of b}
b'(a) = - \frac{\frac{\partial}{\partial a} G(a,b(a))}{\frac{\partial}{\partial b} G(a, b(a))} > 0 
\quad \Leftrightarrow \quad
\frac{\partial}{\partial a} G (a,b(a))>0 ,
\end{align}
and we in order to obtain the desired results, we distinguish two cases.

{\it Case} (a). 
If $a>m$, then we notice that 
\begin{align}\label{Equ: G partial a}
\frac{\partial}{\partial a} G(a, b(a)) = \frac{(\theta_1 -1)(1 - \theta_2)(k(\lambda - (r-g))-\alpha)}{(\theta_1-\theta_2)(\lambda - (r-g))a} \Big[ \Big(\frac{b}{a} \Big)^{1-\theta_1} - \Big(\frac{b}{a}\Big)^{1-\theta_2} \Big] > 0. 
\end{align} 
This implies, thanks to \eqref{Def: Derivative of b}, that 
\begin{equation} \label{bincm+}
a \mapsto b(a) \quad \text{is increasing on} \quad (m,\infty) .
\end{equation}
Furthermore, combining the expressions of the partial derivatives in \eqref{Equ: G partial b} and \eqref{Equ: G partial a} with the expression of $b'(\cdot)$ in \eqref{Def: Derivative of b} yields  $b'(a) = \frac{b(a)}{a}$, which further implies that $b(a) = (1/\tilde{q}) a$, for some $\tilde{q} \in (0,1)$. The latter can be specified as the unique equation to the solution $G(\tilde{q},1)=0$, which is equivalent to \eqref{Def: tilde q}.

{\it Case} (b). 
If $a \leq m$, then we notice that 
\begin{align*}
\frac{\partial}{\partial a} &G(a,b(a)) \\
&= \frac{2 (\alpha - \kappa (\lambda - (r-g))}{(\theta_1 - \theta_2) a \sigma^2} \Big[ \Big(\frac{b(a)}{a}\Big)^{1- \theta_2} - \Big(\frac{b(a)}{a}\Big)^{1- \theta_1} \Big] + \frac{2 \alpha}{(\theta_1 - \theta_2) \sigma^2 a} \Big[ \Big( \frac{m}{a}\Big)^{1 - \theta_1} - \Big(\frac{m}{a}\Big)^{ 1 - \theta_2} \Big] \\
&\geq \frac{2 (\alpha - \kappa (\lambda - (r-g))}{(\theta_1 - \theta_2) a \sigma^2} \Big[ \Big(\frac{b_0}{a}\Big)^{1- \theta_2} - \Big(\frac{b_0}{a}\Big)^{1- \theta_1}\Big] + \frac{2 \alpha}{(\theta_1 - \theta_2) \sigma^2 a} \Big[ \Big( \frac{m}{a}\Big)^{1 - \theta_1} - \Big(\frac{m}{a}\Big)^{ 1 - \theta_2} \Big] \\
&= \Big[ \frac{2 \alpha}{(\theta_1-\theta_2)\sigma^2 a}  - \frac{2 (\alpha - \kappa ( \lambda - (r-g)) b_0^{1- \theta_1} }{(\theta_1-\theta_2) \sigma^2 a} \Big] \geq 0,
\end{align*}
where the first inequality follows from $b(a) \geq b_0$ and the second one from its definition \eqref{Def: b0}.
This implies, thanks to \eqref{Def: Derivative of b}, that 
\begin{equation} \label{bincm-}
a \mapsto b(a) \quad \text{is increasing on} \quad (0, m].
\end{equation}

\vspace{1mm}
\textit{Proof of part} (iv). Regarding the convexity of $b(\cdot)$ on the interval $(0,m)$, we examine the term 
\begin{align*}
b''(a) = \frac{2 G_a (a,b(a)) G_b(a,b(a)) G_{ab}(a,b(a)) - G_b^2(a,b(a)) G_{aa}(a,b(a)) - G_a^2(a,b(a)) G_{bb}(a,b(a))}{G_b^3(a,b(a))},
\end{align*}
where $G_a (a,b) := \partial_a G(a,b)$ (and the other terms analogously). 
We first notice that, due to \eqref{Equ:GPartb}, we have $G_b(a, b(a))<0$. 
Upon using \eqref{Def: G(a,b)} and some direct calculation, we find 
\begin{align*} %\label{Equ: Numerator b''}
&2 G_a (a,b(a)) G_b(a,b(a)) G_{ab}(a,b(a)) - G_b^2(a,b(a)) G_{aa}(a,b(a)) - G_a^2(a,b(a)) G_{bb}(a,b(a)) \nonumber \\
&=\frac{(\theta_1-1)^3 (1-\theta_2)^3 (\alpha - \kappa(\lambda - (r-g))\alpha}{(\lambda - (r-g))^3 (\theta_1-\theta_2)^3 a^2b(a)^2} \nonumber \\
&\quad \times \Big\{ (\alpha - \kappa(\lambda - (r-g)) \Big[ \Big( \frac{b(a)}{a}\Big)^{1-\theta_1} - \Big(\frac{b(a)}{a}\Big)^{1- \theta_2} \Big]^2 \Big[ \theta_2 \Big( \frac{m}{a}\Big)^{1- \theta_2} - \theta_1 \Big( \frac{m}{a}\Big)^{1-\theta_1} \Big] \nonumber \\ 
&\quad \qquad+ \alpha \Big[ \Big( \frac{m}{a}\Big)^{1-\theta_2} - \Big( \frac{m}{a} \Big)^{1-\theta_1} \Big]^2 \Big[ \theta_1 \Big(\frac{b(a)}{a}\Big)^{1-\theta_1} - \theta_2 \Big(\frac{b(a)}{a}\Big)^{1-\theta_2} \Big] \Big\} \nonumber
\end{align*}
While the first term on the above right-hand side (second line) %of \eqref{Equ: Numerator b''} 
is clearly positive, one can employ the fact that $b(a)>b_0>m$ for all $a>0$, with $b_0$ as in \eqref{Def: b0}, to show that the second term is strictly negative. 
Consequently, we obtain $b''(a) > 0$ and thus the strict convexity of $b(\cdot)$ on $(0,m)$.
    
Moreover, we straightforwardly calculate $\lim_{a \downarrow 0} G(a,b_0) = 0$, which due to the the reasoning above and the monotonicity of $b(a)$ implies $\lim_{a \downarrow 0} b(a) = b_0$. \qed

\subsection{Proof of Theorem \ref{VT:Icase2}}\label{Appendix:ProofVTI2}
We derive the result in a number of steps. 
In particular, we show in {\it Step 1} that the candidate value function $U_2(x;a)$ of \eqref{Def: Candidate VF U2} solves the free-boundary problem \eqref{Def: Free Boundary problem ECB}, with $b(a)$ solving \eqref{Equ: G=0} as in Lemma \ref{Lem:b(a)}, and satisfies the HJB equation \eqref{Def: HJB ECB}, and in {\it Steps 2--4} that the latter identifies with value function $V_2$ of \eqref{Value function V2 (x,a)}.  

\vspace{1mm}
{\it Step 1.}
By construction, we have $(\mathcal{L} - \lambda) U_2 (x;a) + \alpha (x-m)^+ =0$ for $x \in (a,b(a))$ as well as $\frac{\partial}{\partial x} U_2 (x;a) = \kappa$ for $x \geq b(a)$. 
It thus remains to show that: 
$$
(i) \;(\mathcal{L}- \lambda) U_2 (x;a) + \alpha(x-m)^+ \geq 0 \quad \text{for $x \geq b(a)$} 
\qquad \text{and} \qquad 
(ii) \;U_2'(x;a) \leq \kappa \text{for $x \in (a,b(a))$.}
$$
In the following, we fix $a \in \mathbb{R}_+$, so that $b(a)$ is the (fixed) unique solution to \eqref{Equ: G=0} according to Lemma \ref{Lem:b(a)}.

{\it Proof of (i)}.
For $x \geq b(a)$, we notice from the expression \eqref{Def: Candidate VF U2} of $U_2$ that  
\begin{align*}
(\mathcal{L}- \lambda) U_2(x;a) 
= (\mathcal{L}- \lambda)\Big( U_2(b(a);a) + \kappa (x-b(a)) \Big) 
= (r-g)x \kappa - \lambda U_2 (b(a);a) - \lambda \kappa (x-b(a)) 
\end{align*}
which is a decreasing, linear function in $x$, with slope $- \kappa (\lambda - (r-g))$. 
We then observe that $x \mapsto - \alpha (x-m)$ also decreases linearly for $x \geq b(a) > m$, with a slope $- \alpha$ that is considered to satisfy $-\alpha < -\kappa (\lambda -(r-g))$, according to the parameter regime under Case (II).
The claim thus follows by noticing that for $x=b(a)$, we have $(\mathcal{L} - \lambda)U_2 (b(a);a) = - \alpha (b(a) - m)^+ =  - \alpha (b(a) - m) $. 

{\it Proof of (ii)}.
For $x \in (a,b(a))$, we notice that $U_2'(x;a) = G(x, b(a))$. 
Since $x \mapsto G(x,b(a))$ is increasing (recall \eqref{Def: Derivative of b} and Lemma \ref{Lem:b(a)}), we obtain via the representation \eqref{Def: G(a,b)} of $G$ that  
\begin{align*}
U_2'(x;a) = G (x,b(a)) \leq G (b(a),b(a)) = \kappa,
\end{align*}
which concludes our claim. Furthermore, we have $ \frac{\partial}{\partial x} U_2 (x;a) \geq G (a,b(a)) = 0$ for $x >a$. 

\vspace{1mm}
{\it Step 2.} Let $x \geq 0$ and $\eta \in \mathcal{A}_{\xi^a}$ such that $(\xi^a , \eta ) \in \mathcal{A}$. Since $U_2 \in C^2 (a, \infty)$, we can apply It\^{o}-Meyer formula, up to a localising sequence of stopping times given by $\tau_n := \inf\{ t \geq 0: X_t^{\xi^a,0} \geq n \}~\P$-a.s., to the process $U_2 (X^{\xi^a,\eta};a)$ and obtain 
\begin{align} \label{Equ: Ito on U2}
e^{-\lambda \tau_n} U_2(X_{\tau_n}^{\xi^a, \eta} ; a) - U_2(x;a) 
= &\int_0^{\tau_n} (\mathcal{L} - \lambda) U_2 (X_s^{\xi^a, \eta} ;a) ds + \sigma \int_0^{\tau_n} U_2'(X_s^{\xi^a, \eta} ;b) dW_s \nonumber \\
& - \int_0^{\tau_n} e^{- \lambda s} X_s^{\xi^a, \eta} U_2'(X_s^{\xi^a, \eta} ;a) d \eta_s^c + \int_0^{\tau_n} e^{- \lambda s} X_s^{\xi^a, \eta} U_2'(X_s^{\xi^a, \eta};a)d \xi_s^{c,a} \nonumber \\
&+ \sum_{s < \tau_n} e^{- \lambda s} \big(U_2 (X_{s}^{\xi^a, \eta} ; a) - U_2 (X_{s-}^{\xi^a, \eta} ; a) \big) 
\end{align}
The second term in \eqref{Equ: Ito on U2} is a martingale due to the continuity of $U_2'$ and the fact that $a \leq X_s^{\xi^a, \eta} < n$ for any $s \in (0, \tau_n]$. Furthermore, we can proceed similarly as in \eqref{Equ: Sum for circ u notation} %-\eqref{Equ: Sum for circ d notation} 
in order to rewrite the last term in \eqref{Equ: Ito on U2} and obtain, after taking expectations and rearranging terms, that
\begin{align}\label{Equ: Inequality VT b}
U_2 (x;a) &= \E_x \Big[ e^{- \lambda \tau_n} U_2 (X_{\tau_n}^{\xi^a, \eta} ;a) - \int_0^{\tau_n} e^{- \lambda s} (\mathcal{L} - \lambda)U_2(X_s^{\xi^a, \eta} ; a) ds \nonumber \\
&\hspace*{1cm}+ \int_0^{\tau_n} e^{- \lambda s} X_s^{\xi^a, \eta} U_2' (X_s^{\xi^a, \eta} ;a) \circ_d d \eta_s - \int_0^{\tau_n} e^{- \lambda s} X_s^{\xi^a, \eta} U_2'(X_s^{\xi^a, \eta} ;a) \circ_u d \xi_s^a \Big] \nonumber \\
&\leq \E_x \Big[  e^{- \lambda \tau_n} U_2 (X_{\tau_n}^{\xi^a, \eta} ;a) + \int_0^{\tau_n} e^{- \lambda s} \alpha (X_s^{\xi^a, \eta} - m)^+ ds + \kappa\int_0^{\tau_n} e^{- \lambda s} X_s^{\xi^a, \eta}  \circ_d d \eta_s   \Big], 
\end{align}
where the latter inequality follows from the fact that $U_2 $ solves the free-boundary problem \eqref{Def: Free Boundary problem ECB} and $U_2 ' (X_s^{\xi^a, \eta};a ) = 0$ for all $s$ in the support of $d \xi^a$. By admissibility of $\eta$, the right-hand side of \eqref{Equ: Inequality VT b} is finite $\P$-a.s., and Assumption \ref{Assumption: c1 > c2}.(iii) %{Assumption: lambda} guarantees 
\begin{align*}
\lim_{n \to \infty} \E_x \big[ e^{- \lambda \tau_n} U_2 (X_{\tau_n}^{\xi^a,\eta};a) \big] = 0.
\end{align*}
Then, taking limits in \eqref{Equ: Inequality VT b} upon using that $\tau_n \uparrow \infty$, we can employ dominated convergence due to Assumption \ref{Assumption: c1 > c2}.(iii) %\ref{Assumption: lambda} 
and obtain  
\begin{align*}
U_2 (x;a) \leq \E_x \Big[  \int_0^{\infty} e^{- \lambda s} \alpha (X_s^{\xi^a, \eta} - m)^+ ds + \kappa\int_0^{\infty} e^{- \lambda s} X_s^{\xi^a, \eta}  \circ_d d \eta_s   \Big].
\end{align*}
We conclude that $U_2 (x;a) \leq V_2 (x;a)$ on $\mathbb{R}_+$.  

\vspace{1mm}
{\it Step 3.} We can now repeat the arguments from {\it Step 2}, upon fixing the control strategy $\eta^{b(a)}$ of \eqref{Def: Control eta b(a)}. Since $X_t^{\xi^a, \eta^{b(a)}} \in [a, b(a)]$ a.s.~for all $t > 0$ and $U_2(X_t^{\xi^a, \eta^{b(a)}};a) = \kappa$ on supp$\{d \eta^{b(a)} \}$, the inequality in \eqref{Equ: Inequality VT b} becomes an equality. 
Arguing as before, we thus obtain  
\begin{align*}
U_2 (x;a) = \E_x \Big[ \int_0^\infty e^{- \lambda s} \alpha (X_s^{\xi^a,\eta^{b(a)}} - m)^+ ds + \kappa \int_0^\infty e^{- \lambda s} X_s^{\xi^a,\eta^{b(a)}} \circ_d d \eta^{b(a)}_s \Big] \geq V_2 (x;a) .
\end{align*}

\vspace{1mm}
{\it Step 4.} Combining the results from {\it Steps 2} and {\it 3} then concludes that $U_2 (x;a) = V_2(x;a)$ on $\mathbb{R}_+$ and $\eta^{b(a)}$ is an optimal control strategy in problem \eqref{Value function V2 (x,a)}.

\section{Proofs of results in Section \ref{nash}}
\subsection{Proof of Theorem \ref{Theorem: Existence and Uniqueness of NE}}\label{Appendix:NashEx}
We prove separately the existence and uniqueness of the Nash Equilibrium. 

\subsubsection*{Existence of a Nash Equilibrium} 
\begin{figure}[h]
    \centering
    \includegraphics[height=4cm]{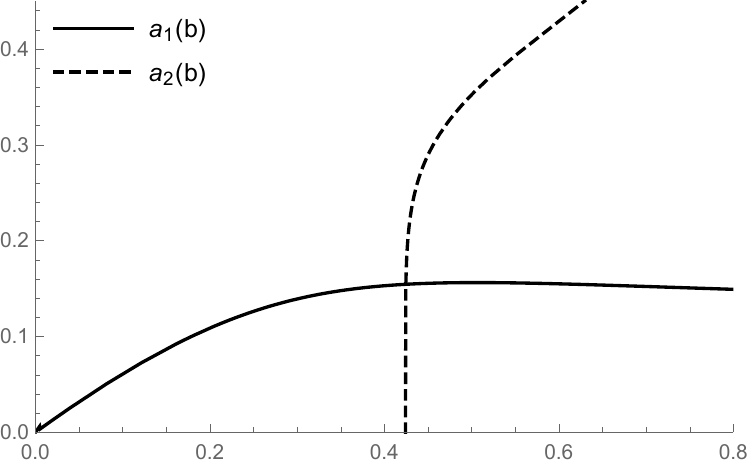}~~
    \includegraphics[height=4cm]{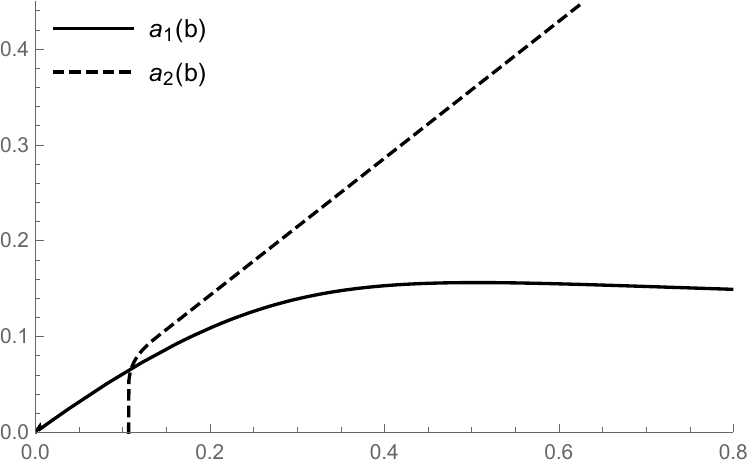}
    \caption{A Sketch of the maps $a_1(b)$ and $a_2(b)$ for different parameter specifications.}
    \label{fig: a1 and a2}
\end{figure}

Recall the function $a(b)$ from Lemma \ref{Lemma: a(b) unique} and define the function \begin{equation} \label{a1} 
a_1(x) := a(x), \quad \text{for all $x \in \mathbb{R}_+$}.
%a_1(b) := a(b), \quad \text{for all $b \in \mathbb{R}_+$}.
\end{equation}
Then, we recall from Lemma \ref{Lem:b(a)} that the unique solution to $G(a,\cdot) = 0$ for any fixed $a \in \mathbb{R}_+$, is given in terms of a strictly increasing function $a \mapsto b(a)$. 
We can therefore invert this function and define 
\begin{equation} \label{a2} 
a_2(x) := b^{-1}(x), \quad \text{such that $x \mapsto a_2(x)$ is strictly increasing on $\mathbb{R}_+$.}
%a_2(b) := \nr{b^{-1}(b)}, %b^{-1}(a), 
%\quad \text{such that $b \mapsto a_2(b)$ is strictly increasing on $\mathbb{R}_+$.} 
\end{equation}
Thanks to Proposition \ref{Lem:b(a)}, we also know that $b(a) = (1/\tilde{q}) a$, for all $a >m$, which yields that 
\begin{equation} \label{a2linear} 
a_2(x) = \tilde{q} x, \quad \text{ for $x > m/\tilde{q}$}.
\end{equation}
For illustration, Figure \ref{fig: a1 and a2} sketches the maps for different parameter specifications.
We can thus conclude that 
\begin{equation} \label{Nasha1a2} 
\exists \text{ an intersection point $b^*$, such that $a^*:=a_1(b^*)=a_2(b^*)$} 
\quad \Rightarrow \quad 
\exists \text{ a Nash equilibrium.} 
\end{equation}
In view of the definitions in \eqref{a1}--\eqref{a2}, we have that $a^*=a(b^*)$ and $b^*=b(a^*)$.
This would finally imply that $(a^*,b^*)$ solves the system of equations $F(a^*,b^*) = 0 = G(a^*,b^*)$ and complete the proof.

In the remainder of the existence proof, we show that the $b^*$ in \eqref{Nasha1a2} indeed exists. 
On one hand, it follows from \eqref{a1} and \eqref{a'0-inf} in Lemma \ref{Lemma: a(b) unique} that $a_1(\cdot)$ is bounded from above by $\hat{a}:=a(\hat{b})$.
On the other hand, it follows from \eqref{a2}--\eqref{a2linear} that $a_2(\cdot)$ is strictly increasing with $\lim_{b \to \infty} a_2(b) = + \infty$. 
We further know from Lemmas \ref{Lemma: a(b) unique} and \ref{Lem:b(a)} that the functions $a_1(\cdot)$ and $a_2(\cdot)$ have supports on $(0,\infty)$ and $(b_0, \infty)$, respectively, with $b_0 > m > 0$. due to \eqref{Def: b0}.
Clearly, there exists at least one $b^* \in (b_0, \infty)$ such that $a_1(b^*) = a_2(b^*)$, therefore \eqref{Nasha1a2} implies that there exists a Nash equilibrium. 

\subsubsection*{Uniqueness of Nash Equilibrium in the class $\mathcal{M}$}
   % Recall from Theorem \ref{Theorem: Existence and Uniqueness of NE} that there exists an intersection point $b^*$, such that $a^*:=a_1(b^*)=a_2(b^*)$, thus $(a^*,b^*)$ is a Nash equilibrium (cf.~\eqref{Nasha1a2}).
In order to prove the uniqueness of the equilibrium established in the previous step, we must prove the uniqueness of the intersection point $b^*$.

We begin by defining the function 
$$
a \mapsto b_0(a), \quad \text{for } a \in \mathbb{R}_+ , 
\quad \text{such that } b_0 (a) = (1/\tilde{q}) a 
\quad \text{with $\tilde{q} \in (0,1)$ as in Lemma \ref{Lem:b(a)}},
$$ 
which can be inverted to define 
$$
a_0 (b) := \tilde{q}b , \quad \text{for } b\in\mathbb{R}_+.
$$
Notice from Lemma \ref{Lem:b(a)} that $b(a) = b_0(a)$ for all $a>m$, hence in view of \eqref{a2}, we get $a_2(b) = a_0(b)$, for all $b>m/\tilde{q}$. Moreover, Lemma \ref{Lem:b(a)} implies that $b \mapsto a_2(b)$ is strictly concave on $(b_0, \frac{m}{\tilde{q}})$, where we notice that $m/\tilde{q} > b_0$ due to the monotonicity of $b(a)$. Moreover, this implies $a_2(b) \leq a_0(b)$ for all $b \geq b_0$. 

As a first step, we prove that the curves $a_1(b)$ and $a_0(b)$ either admit no intersection or exactly one for $b>0$. Since $a_0(b) = \tilde{q} b$, any intersection clearly is of the form $(\tilde{q}b,b)$. Plugging in points of this form into the function $F$ of \eqref{Def: F(a,b)} we observe that
\begin{align*}
b \to F(\tilde{q}b,b) = 
s(\tilde{q}) b + y(\tilde{q}) 
\quad \text{is strictly increasing on $\mathbb{R}_+$},
\end{align*}
with
\begin{align*}
s(q) &= (2 - \delta_2) {q}^{2-\delta_1} + (\delta_1-2) {q}^{2-\delta_2} - (\delta_1 - \delta_2) > 0, \\
y(q) &= (\rho - 2(r-g)-\sigma^2) \big[ c_1 (\delta_1-\delta_2) - c_2 (1-\delta_2)q^{1-\delta_1} - c_2(\delta_1-1)q^{1-\delta_2} \big] .
\end{align*} 
Clearly, this implies (depending on the value $y(q)$) that either no intersection point exists or exactly one. Moreover, if no intersection of $a_0$ and $a_1$ exists, we conclude that $a_1'(0+)< a_0'(0+) = \tilde{q}$.

Next, we come back to the uniqueness of an intersection of the maps $a_1(b)$ and $a_2(b)$. Recall that $b \mapsto a_1(b)$ is concave on $(0,\hat{b})$, and $a \mapsto b(a)$ is convex, which implies that $b \mapsto a_2(b)$ is concave as well.  

We now distinguish the following cases; an illustration for each one of them is given in Figure \ref{Fig:CasesUni}:

\begin{figure}[h]
\centering
\begin{subfigure}{0.35\textwidth}
\includegraphics[height=3.5cm]{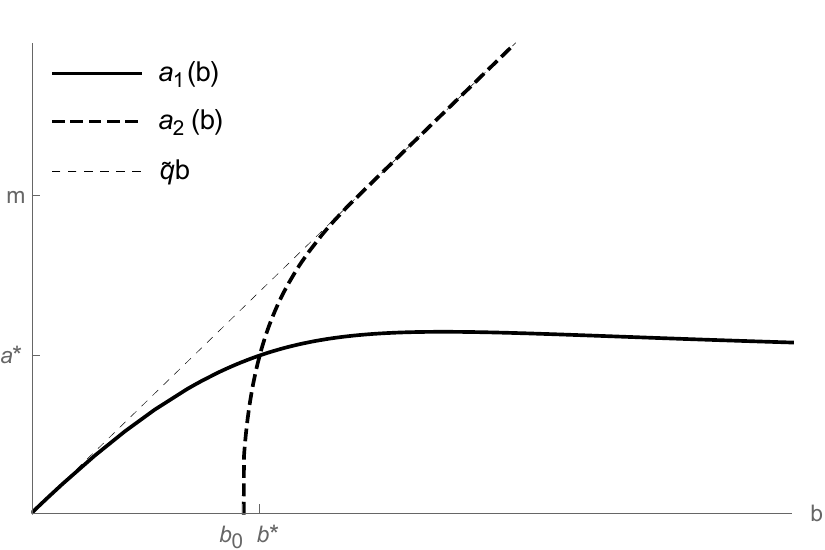} \caption*{Case (i)}
\label{fig:subcase(i)}
\end{subfigure}\begin{subfigure}{0.35\textwidth}
\includegraphics[height=3.5cm]{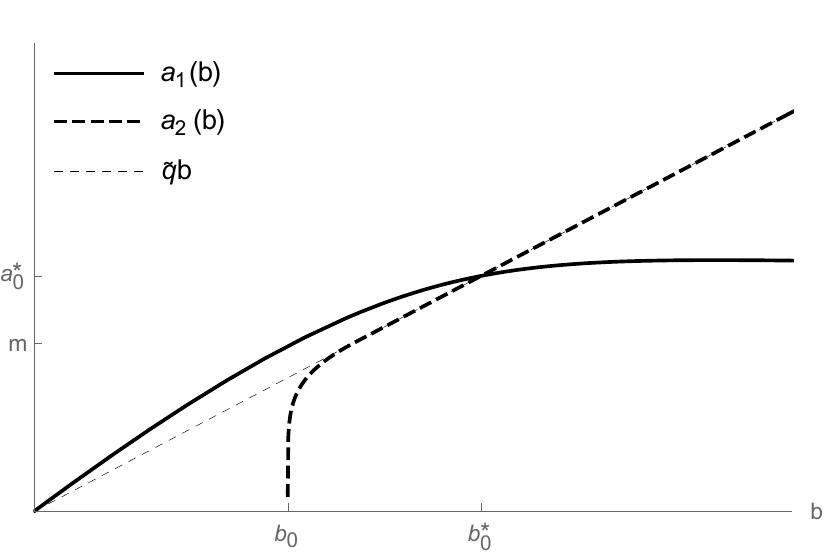}
    \caption*{Case (ii)}
\label{fig:subcase(ii)}
\end{subfigure}

\begin{subfigure}{0.35\textwidth}
\includegraphics[height=3.5cm]{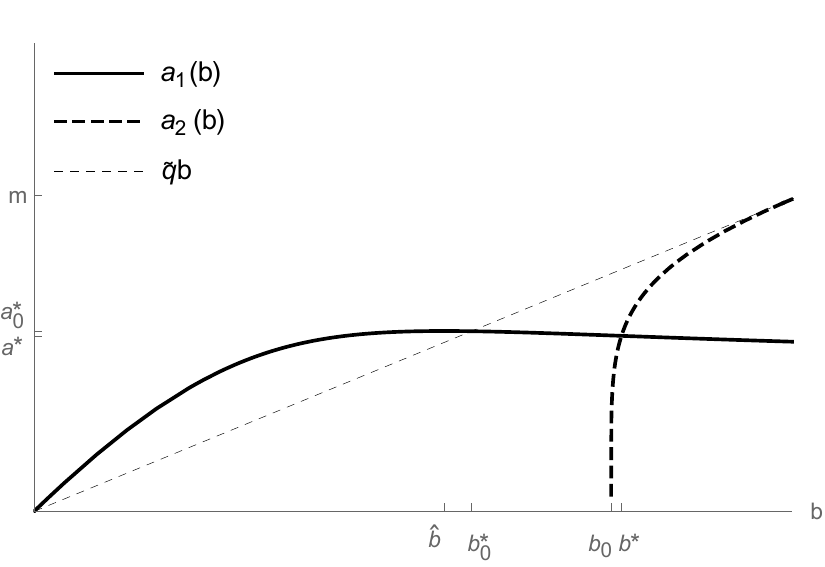}
    \caption*{Case (iii)}
\label{fig:subcase(iii)}
\end{subfigure}\begin{subfigure}{0.35\textwidth}
\includegraphics[height=3.5cm]{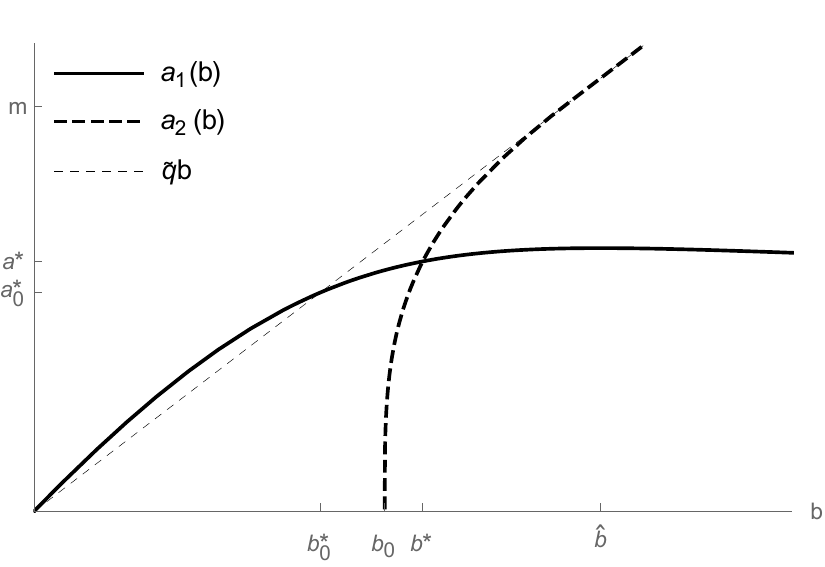}
    \caption*{Case (iv)}
\label{fig:subcase(iv)}
\end{subfigure}
\caption{Case Study for the uniqueness of an intersection of $a_1(b)$ and $a_2(b)$}
\label{Fig:CasesUni}
\end{figure}

{\it Case} (i). No intersection exists of $a_1(b)$ and $a_0(b)$: In this case, we first notice that $a_2'(b) \geq a_0'(b) = \tilde{q}$ for all $b$ in the support of $a_2(b)$ due to the concavity of $a_2(b)$ and the fact that $a_2(b) \leq \tilde{q}b$. Furthermore, we have $\tilde{q} > a_1'(0+)$, and hence $\tilde{q} > a_1'(b)$ for all $b>0$. Combining these insights, it follows that $a_1(b) = a_2(b)$ for exactly one $b\in \mathbb{R}_+$.  

For the following cases, we thus assume that there exists exactly one intersection of $a_1(b)$ and $a_0(b)$, which yields the point $(a^*_0,b^*_0)$. 

{\it Case} (ii). If $a^*_0 \geq m$, the uniqueness of an intersection of $a_1(b)$ and $a_2(b)$ follows from the fact that $a_2(b) = a_0(b)$ for all $b \geq m/\tilde{q}$. Hence, in this case, we obtain that the intersection of $a_1(b)$ and $a_2(b)$ is exactly given by $(a^*, b^*) = (a^*_0,b^*_0)$.

{\it Case} (iii). If $a^*_0 < m$ and $b^*_0 \geq \hat{b}$, we again observe that $a_2(b) \leq \tilde{q}b$, which implies that the intersection point is not realized on $b \leq b_0^*$. Therefore, since $a_1(b)$ is decreasing for $b \geq \hat{b}$ and $a_2(b)$ is increasing, the intersection is unique and we denote it by $(a^*,b^*)$.

{\it Case} (iv). If $a^*_0 < m$ and $b_0^* < \hat{b}$, the concavity of $a_1(b)$ implies $a_1'(b_0^*) < \tilde{q}$ as well as $a_1'(b_0^*) > a_1'(b)$ for all $b \geq b_0^*$. Since $a_2'(b) > \tilde{q}$ for all $b$, we must have $a_1(b) = a_2(b)$ for exactly one $b \in \mathbb{R}_+$ and we again denote the unique intersection by $(a^*,b^*)$. \qed

\subsection{Proof of Corollary \ref{Coro:ConvNash}}\label{Appendix:Coro}
Recall that $a_2(b; \lambda) \leq a_0(b; \lambda) = \tilde{q}(\lambda) b$ (see also Figure \ref{Fig:CasesUni} for illustration), where $\Tilde{q} \equiv \tilde{q}(\lambda) \in (0,1)$ is given by the solution to \eqref{Def: tilde q}, which is equivalent to $G(\Tilde{q}(\lambda),1)=0$. 
By defining %Let 
\begin{align*}
\hat{G}(q) := 
(1 - \theta_2) \kappa (\lambda - (r-g))   + (\theta_1 -1) (\kappa (\lambda - (r-g)) - \alpha) q^{\theta_2 -1} + \alpha (\theta_1 -1), \qquad q \in (0,1),
\end{align*}
%for $q \in (0,1)$, 
we observe that $G(q,1) \geq \hat{G}(q)$ for all $q\in (0,1)$, since $\lambda < r-g + \frac{\alpha}{\kappa}$. 
Moreover, we notice that 
\begin{align*}
\lim_{q \downarrow 0} G(q,1) = - \infty, \quad \lim_{q \downarrow 0} \hat{G}(q) = - \infty , \quad G(1,1) = \hat{G}(1) = (\theta_1 - \theta_2)\kappa (\lambda - (r-g)) > 0
\end{align*} 
and that both $q \mapsto G(q,1)$ and $q \mapsto \hat{G}(q)$ are monotonically increasing. 
Using these properties, we can denote by $\hat{q} \equiv \hat{q}(\lambda) \in (0,1)$ the unique solution to $\hat{G}(q) = 0$ (which can be computed explicitly) and 
obtain  
\begin{align*} 
\tilde{q}(\lambda) < \hat{q}(\lambda) \quad \text{and} \quad 
\lim_{\lambda \uparrow r-g + \alpha/\kappa} \hat{q}(\lambda) = 0 
\quad \Rightarrow \quad 
\lim_{\lambda \uparrow r-g + \alpha/\kappa} \tilde{q}(\lambda) = 0.
%\quad \Rightarrow \quad \lim_{\lambda \uparrow r-g + \alpha/\kappa} a_2(b; \lambda) = 0
\end{align*}

First of all, recall from Theorem \ref{Theorem: Existence and Uniqueness of NE} that for every $\lambda \in (r-g, r-g+\alpha/\kappa)$, there exists a unique equilibrium pair $(a^*(\lambda),b^*(\lambda))$, such that $a^*(\lambda) = a_1(b^*(\lambda)) = a_2(b^*(\lambda);\lambda)$.  
Since $a_2(b;\lambda) \leq \tilde{q}(\lambda) b$, we notice that for any fixed $\Tilde{b} \in (0,\infty)$ we have $\lim_{\lambda \uparrow r-g+\alpha/\kappa} a_2(\tilde{b};\lambda) =0$. 
However, there is no such fixed $\Tilde{b} \in (0,\infty)$ that can give $a_1(\Tilde{b}) = \lim_{\lambda \uparrow r-g+\alpha/\kappa} a_2(\Tilde{b};\lambda) = 0$ to create an equilibrium pair $\lim_{\lambda \uparrow r-g+\alpha/\kappa} (a^*(\lambda), b^*(\lambda)) = (0, \Tilde{b})$. 
Hence, the only possibility for obtaining an equilibrium as $\lambda \uparrow r-g+\alpha/\kappa$, is for $b^*(\lambda) \to \infty$.
Given that $b \mapsto a_1(b)$ is strictly decreasing on $(\hat{b},\infty)$ and $\lim_{b \to \infty} a_1(b) = \overline{a}$ (independently of $\lambda$), we conclude that 
$$
\lim_{\lambda \uparrow r-g+\alpha/\kappa} a^*(\lambda) = \lim_{\lambda \uparrow r-g+\alpha/\kappa} a_2(b^*(\lambda);\lambda) = \lim_{\lambda \uparrow r-g+\alpha/\kappa} a_1(b^*(\lambda)) = \overline{a}
\quad \text{and} \quad 
\lim_{\lambda \uparrow r-g+\alpha/\kappa} b^*(\lambda) = \infty,
$$ 
which completes the proof. \qed

\section*{Declarations}
\text{}\\[-0.3cm]\textbf{Conflict of interest} The authors have no conflicts of interest to declare that are relevant to the content of
this article.

\end{document}